\documentclass{gtmon_a}
\pdfoutput=1

\usepackage{amscd}
\usepackage{subfigure}
\usepackage[arrow,matrix,graph,frame,poly,arc,tips]{xy}

\proceedingstitle{Groups, homotopy and configuration spaces (Tokyo
  2005)}
\conferencestart{5 July 2005}
\conferenceend{11 July 2005}
\conferencename{Groups, homotopy and configuration spaces, 
                in honour of Fred Cohen's 60th birthday}
\conferencelocation{University of Tokyo, Japan}

\editor{Norio Iwase}
\givenname{Norio}
\surname{Iwase}

\editor{Toshitake Kohno}
\givenname{Toshitake}
\surname{Kohno}

\editor{Ran Levi}
\givenname{Ran}
\surname{Levi}

\editor{Dai Tamaki}
\givenname{Dai}
\surname{Tamaki}

\editor{Jie Wu}
\givenname{Jie}
\surname{Wu}

\title{The boundary manifold of a complex line arrangement}

\author{Daniel C Cohen}
\givenname{Daniel C}
\surname{Cohen}
\address{Department of Mathematics\\
Louisiana State University\\\newline
Baton Rouge LA 70803\\
USA}
\email{cohen@math.lsu.edu}
\urladdr{http://www.math.lsu.edu/~cohen/}

\author{Alexander I Suciu}
\givenname{Alexander I}
\surname{Suciu}
\address{Department of Mathematics\\
Northeastern University\\\newline
Boston MA 02115\\
USA}
\email{a.suciu@neu.edu}
\urladdr{http://www.math.neu.edu/~suciu/}

\dedicatory{For Fred Cohen on the occasion of his sixtieth birthday}

\volumenumber{13}
\issuenumber{}
\publicationyear{2008}
\papernumber{5}
\startpage{105}
\endpage{146}

\doi{}
\MR{}
\Zbl{}

\arxivreference{math.GT/0607274}


\received{30 May 2006}
\revised{29 May 2007}
\accepted{30 May 2007}
\published{22 March 2008}
\publishedonline{22 February 2008}
\proposed{}
\seconded{}
\corresponding{}
\version{}

\makeatletter
\def\cnewtheorem#1[#2]#3{\newtheorem{#1}{#3}[section]
\expandafter\let\csname c@#1\endcsname\c@subsection}

\let\xysavmatrix\xymatrix
\def\xymatrix{\disablesubscriptcorrection\xysavmatrix}
\let\xysavgraph\xygraph
\def\xygraph{\disablesubscriptcorrection\xysavgraph}
\AtBeginDocument{\let\bar\wbar\let\tilde\wtilde\let\hat\what}

\theoremstyle{plain}
\cnewtheorem{thm}[subsection]{Theorem}
\cnewtheorem{prop}[subsection]{Proposition}
\cnewtheorem{lem}[subsection]{Lemma}
\cnewtheorem{cor}[subsection]{Corollary}
\cnewtheorem{conj}[subsection]{Conjecture}

\theoremstyle{definition}
\cnewtheorem{definition}[subsection]{Definition}
\cnewtheorem{example}[subsection]{Example}

\theoremstyle{remark}
\newtheorem*{remark}{Remark}
\newtheorem*{rem}{Added in proof}
\newtheorem*{ack}{Acknowledgment}

\makeatother  %  move after \newtheorem block

\newcommand{\abs}[1]{\left|#1\right|}
\newcommand{\set}[1]{\left\{#1\right\}}

\renewcommand{\b}[1]{\mathbf{#1}}   
\newcommand{\norm}[1]{\|#1\|}

\renewcommand{\atop}[2]{\genfrac{}{}{0pt}{}{#1}{#2}}

\newcommand{\surj}{{\twoheadrightarrow}}

\newcommand{\bB}{\mathbb{B}}
\newcommand{\bF}{\mathbb{F}}
\newcommand{\bS}{\mathbb{S}}
\newcommand{\F}{\mathbb{F}}
\newcommand{\CP}{{\mathbb{CP}}}

\newcommand{\sD}{{\sf D}}
\newcommand{\GA}{{\Gamma_{\!\!\A}}}

\newcommand{\A}{\mathcal{A}}
\newcommand{\RR}{\mathcal{R}}

\newcommand{\cV}{\mathcal{V}}
\newcommand{\cE}{\mathcal{E}}
\newcommand{\cC}{\mathcal{C}}
\newcommand{\cT}{\mathcal{T}}
\newcommand{\cF}{\mathcal{F}}
\newcommand{\cN}{\mathcal{N}}
\newcommand{\cS}{\mathcal{S}}

\DeclareMathOperator{\rank}{rank}

\DeclareMathOperator{\Hom}{Hom}
\DeclareMathOperator{\Image}{Im}
\DeclareMathOperator{\GL}{GL}
\DeclareMathOperator{\Sym}{Sym}

\DeclareMathOperator{\ord}{ord}
\DeclareMathOperator{\Tors}{Tors}
\DeclareMathOperator{\Int}{Int}

\DeclareMathOperator{\gr}{gr}
\DeclareMathOperator{\ab}{ab}
\DeclareMathOperator{\Hilb}{Hilb}

\DeclareMathOperator{\PD}{PD}
\DeclareMathOperator{\cham}{ch}
\DeclareMathOperator{\Lie}{{Lie}}
\DeclareMathOperator{\nbc}{\mathbf{nbc}}

\newcommand{\h}{{\mathfrak{h}}}
\newcommand{\chim}{{\chi_{\mathunderscore}}}
\newcommand{\db}[1]{{\what{#1}}} 
\newcommand{\dbl}[1]{{\what{#1}}}

\newcommand{\dA}{{\mathsf{d}\mathcal{A}}}
\newcommand{\disc}{\text{\:\:\circle*{0.35}}}
\newcommand{\disca}{\text{\:\circle*{0.35}}}
\newcommand{\discb}{\text{\:\:\:\circle*{0.35}}}

\newenvironment{romenum}
{\renewcommand{\theenumi}{\roman{enumi}}

\begin{enumerate}}{\end{enumerate}}

\begin{document}

\begin{htmlabstract}
We study the topology of the boundary manifold of a line arrangement
in <b>CP</b><sup>2</sup>, with emphasis on the fundamental group G and
associated invariants.  We determine the Alexander polynomial &Delta;(G),
and more generally, the twisted Alexander polynomial associated to
the abelianization of G and an arbitrary complex representation.
We give an explicit description of the unit ball in the Alexander norm,
and use it to analyze certain Bieri&ndash;Neumann&ndash;Strebel invariants of G.
From the Alexander polynomial, we also obtain a complete description
of the first characteristic variety of G. Comparing this with the
corresponding resonance variety of the cohomology ring of G enables
us to characterize those arrangements for which the boundary manifold
is formal.
\end{htmlabstract}

\begin{abstract}
We study the topology of the boundary manifold of a line arrangement
in $\mathbb{CP}^2$, with emphasis on the fundamental group $G$ and
associated invariants.  We determine the Alexander polynomial $\Delta(G)$,
and more generally, the twisted Alexander polynomial associated to
the abelianization of $G$ and an arbitrary complex representation.
We give an explicit description of the unit ball in the Alexander norm,
and use it to analyze certain Bieri--Neumann--Strebel invariants of $G$.
From the Alexander polynomial, we also obtain a complete description
of the first characteristic variety of $G$. Comparing this with the
corresponding resonance variety of the cohomology ring of $G$ enables
us to characterize those arrangements for which the boundary manifold
is formal.
\end{abstract}
\begin{asciiabstract}
We study the topology of the boundary manifold of a line arrangement
in CP^2, with emphasis on the fundamental group G and associated
invariants.  We determine the Alexander polynomial Delta(G), and more
generally, the twisted Alexander polynomial associated to the
abelianization of G and an arbitrary complex representation.  We give
an explicit description of the unit ball in the Alexander norm, and
use it to analyze certain Bieri-Neumann-Strebel invariants of G.  From
the Alexander polynomial, we also obtain a complete description of the
first characteristic variety of G. Comparing this with the
corresponding resonance variety of the cohomology ring of G enables us
to characterize those arrangements for which the boundary manifold is
formal.
\end{asciiabstract}

\maketitle

\section{Introduction}
\label{sec:intro}

\subsection{The boundary manifold} 

Let $\A$ be an arrangement of hyperplanes in the complex 
projective space $\CP^m$, $m>1$.  
Denote by $V= \bigcup_{H\in\A} H$ the corresponding 
hypersurface, and by $X=\CP^m \setminus V$ its complement. 
Among the origins of the topological study of arrangements 
are seminal results of Arnol'd \cite{Ar69} and Cohen \cite{FC}, 
who independently computed the cohomology of the 
configuration space of $n$ ordered points in $\C$, the 
complement of the braid arrangement.  The cohomology 
ring of the complement of an arbitrary arrangement $\A$ 
is by now well known.  It is isomorphic to the Orlik--Solomon 
algebra of $\A$, see Orlik and Terao \cite{OT1} as a general 
reference.

In this paper, we study a related topological space, namely 
the \emph{boundary manifold} of $\A$.  By definition, this 
is the boundary $M=\partial N$ of a regular neighborhood of 
the variety $V$ in $\CP^m$.  Unlike the complement $X$, 
an open manifold with the homotopy type of a CW--complex 
of dimension at most $m$, the boundary manifold $M$ is a 
compact (orientable) manifold of dimension $2m-1$.  

In previous work \cite{CS06}, we have shown that 
the cohomology ring of $M$ is functorially determined 
by that of $X$ and the ambient dimension.  In particular, 
$H_*(M;\Z)$ is torsion-free, and the respective Betti numbers 
 are related by $b_k(M)=b_k(X)+b_{2m-k-1}(X)$. 
So we turn our attention here to another topological invariant, 
the fundamental group.  The inclusion map 
$M \to X$ is an $(m-1)$--equivalence, see Dimca \cite{Dimca}.   
Consequently, 
for an arrangement $\A$ in $\CP^m$ with $m \ge 3$, 
the fundamental group of the boundary is isomorphic 
to that of the complement.  In light of this, we focus 
on arrangements of lines in~$\CP^2$.

\subsection{Fundamental group} 

Let $\A=\{\ell_0, \dots, \ell_n\}$ be a line arrangement in $\CP^2$. 
The boundary manifold $M$ is a graph manifold in the sense 
of Waldhausen \cite{Wa1,Wa2}, modeled on a certain weighted graph 
$\GA$. This structure, which 
we review in \fullref{sec:bdry}, has been used by a number 
of authors to study the manifold $M$.  For instance, 
Jiang and Yau \cite{JY93,JY98} investigate the relationship 
between the topology of $M$ and the combinatorics of $\A$, 
and Hironaka \cite{Hi} analyzes the relationship between 
the fundamental groups of $M$ and $X$.  

If $\A$ is a pencil of lines, then $M$ is a connected sum of 
$n$ copies of $S^1 \times S^2$.  Otherwise, $M$ is aspherical, 
and so the homotopy type of $M$ is encoded in its fundamental group.
Using the graph manifold structure, and a method due 
to Hirzebruch \cite{Hir}, Westlund finds a presentation 
for the group $G=\pi_1(M)$ in \cite{We}.  In \fullref{sec:pi1}, we 
build on this work to find a {\em minimal} presentation 
for the fundamental group, of the form 
\begin{equation}
\label{eq:pi1 intro}
G=\langle x_j, \gamma_{i,k} \mid R_j, R_{i,k}\rangle, 
\end{equation}
where $x_j$ corresponds to a meridian loop around line
$\ell_j$, for $1\le j \le n=b_1(X)$, and $\gamma_{i,k}$ corresponds
to a loop in the graph $\Gamma_{\A}$, indexed by a pair
$(i,k)\in \nbc_2(\dA)$, where $\abs{\nbc_2(\dA)}=b_2(X)$.  The
relators $R_j$, $R_{i,k}$ (indexed in the same way) are certain 
products of commutators in the generators.  In other words, 
$G$ is a commutator-relators group, with both 
generators and relators equal in number to $b_1(M)$.

\subsection{Twisted Alexander polynomial and related invariants}

Since $M$ is a graph manifold, the group $G=\pi_1(M)$ 
may be realized as the fundamental group of a graph of groups.  
In \fullref{sec:AlexPolys} and \fullref{sec:alex poly arr}, this structure 
is used to calculate the twisted Alexander polynomial $\Delta^\phi(G)$ 
associated to $G$ and an arbitrary complex representation 
$\phi\colon G \to \GL_k(\C)$.  In particular, we show that the 
classical multivariable Alexander polynomial, arising from 
the trivial  representation of $G$, is given by
\begin{equation}
\label{eq:delta intro}
\Delta(G) = \prod_{v \in \cV(\GA)} (t_v-1)^{m_v-2},  
\end{equation}
where $\cV(\GA)$ is the vertex set of $\GA$, 
$m_v$ denotes the multiplicity or degree of the vertex $v$, 
and $t_v=\prod_{i\in v} t_i$. 

Twisted Alexander polynomials 
inform on invariants such as the Alexander and Thurston norms, 
and  Bieri--Neumann--Strebel (BNS) invariants.  As such, they are 
a subject  of current interest in $3$--manifold theory.  In the case 
where $G$ is a link group, a number of authors, including 
Dunfield \cite{Dun} and Friedl and Kim \cite{FK}, have used 
twisted Alexander polynomials to distinguish between the 
Thurston and Alexander norms.  This is not possible for 
(complex representations of) the fundamental group of the 
boundary manifold of a line arrangement.  In \fullref{sec:alex balls}, 
we show that the unit balls in the norms on $H^1(G;\R)$ 
corresponding to any two twisted Alexander polynomials 
are equivalent polytopes.  Analysis of the structure of these 
polytopes also enables us to calculate the number of 
components of the BNS invariant of $G$ and the Alexander 
invariant of $G$.

\subsection{Cohomology ring and graded Lie algebras}  

In \fullref{sect:coho}, we revisit the cohomology ring of 
the boundary manifold $M$, in our $3$--dimensional context.  
{F}rom \cite{CS06}, we know that $H^*(M;\Z)$ is isomorphic 
to $\db{A}$, the ``graded double" of $A=H^*(X;\Z)$.  In particular, 
$\db{A}^1=A^1\oplus \bar{A}^2$, where $\bar{A}^k=\Hom(A^k,\Z)$. 
This information allows us to identify the $3$--form 
$\eta_M$ which encodes all the cup-product structure in  
the Poincar\'{e} duality algebra $H^*(M;\Z)$.  If $\{e_j\}$ 
and $\{f_{i,k}\}$ denote the standard bases  
for $A^1$ and $A^2$, then 
\begin{equation}
\label{eq:eta intro}
\eta_M =\sum_{(i,k)\in\nbc_2(\dA)}
e_{I(i,k)} \wedge e_k \wedge \bar{f}_{i,k}, 
\end{equation}
where 
$I(i,k)=\set{j \mid \ell_j \supset \ell_i \cap \ell_k,\ 1\le j \le n}$ 
and $e_J=\sum_{j\in J} e_j$. 

The explicit computations described in \eqref{eq:pi1 intro} 
and \eqref{eq:eta intro} facilitate analysis of two Lie algebras 
attached to our space $M$: the graded Lie algebra $\gr(G)$ 
associated to the lower central series of $G$, and the holonomy 
Lie algebra $\h(\db{A})$ arising from the multiplication map 
$\db{A}^1\otimes \db{A}^1\to \db{A}^2$. For the complement 
$X$, the corresponding Lie algebras are isomorphic 
over the rationals, as shown by Kohno \cite{K}. For the 
boundary manifold, though, such an isomorphism no longer holds, 
as we illustrate by a concrete example in \fullref{sect:formal}.  
This indicates that the manifold $M$, unlike the complement 
$X$, need not be formal, in the sense of Sullivan \cite{Su77}. 

\subsection{Jumping loci and formality}

The non-formality phenomenon identified above is fully 
investigated in \fullref{sect:cjl} and \fullref{sect:formal}  by 
means of two types of varieties attached to $M$:  the 
characteristic varieties $V^1_d(M)$ and the resonance varieties 
$\RR^1_d(M)$.   Our calculation of $\Delta(G)$ recorded in 
\eqref{eq:delta intro} enables us to give a complete description 
of the first characteristic variety of $M$, the set of all characters 
$\phi \in \Hom(G,\C^*)$ for which the corresponding local system 
cohomology group $H^1(M;\C_\phi)$ is non-trivial:
\begin{equation}
\label{eq:v1 intro}
V^1_1(M) = \bigcup_{v \in \cV(\GA),m_v \ge 3} \{t_v-1=0\}.
\end{equation}
The resonance varieties of $M$ are the analogous jumping loci 
for the cohomology ring $H^*(M;\C)$.  Unlike the resonance 
varieties of the complement $X$, the varieties $\RR^1_d(M)$, 
for $d$ sufficiently large, may have non-linear components.  
Nevertheless, the first resonance variety $\RR^1_1(M)$ is 
very simple to describe: with a few exceptions, it is equal to the 
ambient space, $H^1(M;\C)$. Comparing the tangent cone 
to $V^1_1(M)$ at the identity to $\RR^1_1(M)$, and making use 
of a recent result of Dimca, Papadima, and Suciu \cite{DPS}, 
we conclude that the boundary manifold of a line arrangement 
$\A$ is formal precisely when $\A$ is a pencil or a near-pencil. 

\section{Boundary manifolds of line arrangements}
\label{sec:bdry}
Let $\A=\set{\ell_0,\dots,\ell_n}$ be an arrangement of lines 
in $\CP^2$.  The boundary manifold of $\A$ may be realized 
as the boundary of a regular neighborhood of the curve 
$C=\bigcup_{i=0}^n \ell_i$ in $\CP^2$.  In this section, 
we record a number of known results regarding this manifold.

\subsection{The boundary manifold}
\label{subsec:bdry nbhd}

Choose homogeneous coordinates $\b{x}=(x_0 \colon x_1 
\colon x_2)$ on $\CP^2$.  For each $i$, $0\le i \le n$, 
let $f_i = f_i(x_0,x_1,x_2)$ be a linear form which 
vanishes on the line $\ell_i$ of $\A$.  Then $Q=Q(\A)=
\prod_{i=0}^n f_i$ is a homogeneous polynomial of degree $n+1$, 
with zero locus $C$.  The {\em complement} of $\A$ is the open  
manifold $X=X(\A)=\CP^2 \setminus C$.

A closed, regular neighborhood $N$ 
of $C$ may be constructed as follows.  Define 
$\phi\colon\CP^2 \to \R$ by $\phi(\b{x}) = 
\abs{Q(\b{x})}^2\!/\, \norm{\b{x}}^{2(n+1)}$, and 
let $N = \phi^{-1}([0,\delta])$ for $\delta>0$ 
sufficiently small.  Alternatively, triangulate 
$\CP^2$ with $C$ as a subcomplex, and take $N$ 
to be the closed star of $C$ in the second barycentric 
subdivision.  As shown by Durfee \cite{Durfee} in greater 
generality, these approaches yield isotopic neighborhoods, 
independent of the choices made in the respective constructions.
The \emph{boundary manifold}  of $\A$ is the 
boundary of such a regular neighborhood:
\begin{equation}
\label{eq:bdry reg nbhd}
M=M(\A)=\partial N.
\end{equation}
This  compact, connected, orientable $3$--manifold will 
be our main object of study. We start with a couple of simple 
examples.

\begin{example} 
\label{ex:boundary pencil}
Let $\A$ be a pencil of $n+1$ lines in $\CP^{2}$, 
defined by the polynomial $Q=x_1^{n+1}-x_2^{n+1}$. 
The complement $X$ of $\A$ is diffeomorphic to  
$(\C \setminus \set{n\ \text{points}})\times \C$, so 
has the homotopy type of a bouquet of $n$ circles.  
On the other hand, $\CP^{2}\setminus N=
(D^2\setminus \{\text{$n$ disks}\})\times D^{2}$; 
hence $M$ is diffeomorphic to the $n$--fold connected 
sum $\sharp^{n} S^{1}\times S^{2}$.  
\end{example}  

\begin{example} 
\label{ex:boundary near-pencil}
Let $\A$ be a near-pencil of $n+1$ lines in $\CP^{2}$, 
defined by the polynomial $Q=x_0(x_1^n-x_2^n)$. 
In this case, $M=S^1\times\Sigma_{n-1}$, 
where $\Sigma_g=\sharp^{g} S^1\times S^1$ 
denotes the orientable surface of genus $g$,  
see \cite{CS06} and \fullref{example:near-pencil pres}. 
\end{example}  

\subsection{Blowing up dense edges}
\label{subsec:blow up}

A third construction, which sheds light on the structure 
of $M$ as a $3$--manifold, may also be used to obtain the 
topological type of the boundary manifold.  This involves 
blowing up (certain) singular points of $C$.  Before 
describing it, we establish some notation.

An edge of $\A$ is a non-empty intersection of lines of $\A$.  
An edge $F$ is said to be {\em dense} if the subarrangement 
$\A_F=\{\ell_j \in \A \mid F \subseteq \ell_j\}$ 
of lines containing $F$ is not a product arrangement.  
Hence, the dense edges are the lines of $\A$, and the 
intersection points $\ell_{j_1} \cap \ldots \cap \ell_{j_k}$ 
of multiplicity $k \ge 3$.  Denote the set of dense edges 
of $\A$ by $\sD(\A)$, and let $F_1,\dots,F_r$ be the 
$0$--dimensional dense edges.  We will occasionally 
denote the dense edge 
$\bigcap_{j \in J} \ell_j$ by $F_J$.

Blowing up $\CP^2$ at each $0$--dimensional dense edge 
of $\A$, we obtain an arrangement 
$\tilde\A=\{L_i\}_{i=0}^{n+r}$ in $\widetilde{\CP}{}^2$ 
consisting of the proper transforms $L_i=\tilde\ell_i$, 
$0\le i \le n$, of the lines of $\A$, and exceptional 
lines $L_{n+j}=\tilde{F_j}$, $1\le j\le r$, arising 
from the blow-ups.  

By construction, the curve 
$\tilde{C}=\bigcup_{i=0}^{n+r} L_i$ in $\widetilde{\CP}{}^2$ 
is a divisor with normal crossings.  Let $U_i$ be a tubular 
neighborhood of $L_i$ in $\widetilde{\CP}{}^2$.  For 
sufficiently small neighborhoods, we have 
$U_i \cap U_j=\emptyset$ if $L_i \cap L_j=\emptyset$.  
Then, rounding corners, $N(\tilde C) = \bigcup_{i=0}^{n+r} U_i$ 
is a regular neighborhood of $\tilde C$ in $\widetilde{\CP}{}^2$.  
Contracting the exceptional lines of $\tilde\A$ gives rise to 
a homeomorphism $M \cong \partial{N}(\tilde C)$.

\subsection{Graph manifold structure}
\label{subsec:graph manifold}

\begin{figure}%
\subfigure{%
\label{fig:np}%
\begin{minipage}[t]{0.4\textwidth}
\setlength{\unitlength}{18pt}
\begin{picture}(4,5.1)(-3.2,-0.8)
\put(0,0){\line(1,1){4}}
\put(-1,2){\line(1,0){6}}
\put(0,4){\line(1,-1){4}}
\put(0.5,-0,5){\line(0,1){5}}
\put(1.2,-0,5){\makebox(0,0){$\ell_0$}}
\put(4.4,-0.5){\makebox(0,0){$\ell_1$}}
\put(5.5,1.95){\makebox(0,0){$\ell_2$}}
\put(4.4,4.5){\makebox(0,0){$\ell_3$}}
\put(2,2.6){\makebox(0,0){$F$}}
\end{picture}
\end{minipage}
}
\subfigure{%
\label{fig:npgraph}%
\begin{minipage}[t]{0.4\textwidth}
\setlength{\unitlength}{20pt}
\begin{picture}(4,5.1)(-2.9,-2.5)
\xygraph{!{0;<14mm,0mm>:<0mm,14mm>::}
[]*D(3){v_{123}}*-{\blacklozenge}
(
-@{--}[dl]*R(2){v_1}*{\disc}(-@{--}[dr]*U(3){v_0}*{\discb})
,-^(0.6){}[d]*R(2){v_2}*{\discb}(-@{--}[d])
,-^{}[dr]*L(2){v_3}*{\disc}(-@{--}[dl])
)
}
\end{picture}
\end{minipage}
}
\caption{A near-pencil of $4$ lines and 
its associated graph $\Gamma$ (with maximal 
tree $\cT$ in dashed lines)}
\label{fig:nearpencil}
\end{figure}
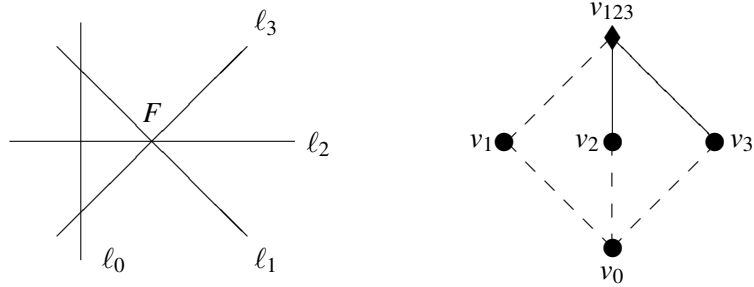

This last construction realizes the boundary manifold $M$ 
of $\A$ as a {\em graph manifold}, in the sense of 
Waldhausen \cite{Wa1,Wa2}.  
The underlying graph $\GA$ may be described as follows.  
The vertex set $\cV(\GA)$ is in one-to-one correspondence 
with the dense edges of $\A$ (that is, the lines of $\tilde\A$).  
Label the vertices of 
$\GA$ by the relevant subsets of $\{0,1,\dots,n\}$:  
the vertex corresponding to $\ell_i$ is labeled $v_i$, 
and, if $F_J$ is a $0$--dimensional dense edge (that is, 
an exceptional line in $\tilde{\A}$), label the 
corresponding vertex $v_J$.  
If $\ell_i$ and $\ell_j$ meet in a double point of $\A$, we 
say that $\ell_i$ and $\ell_j$ are transverse, and (sometimes) 
write $\ell_i\pitchfork\ell_j$.  
The graph $\GA$ has an 
edge $e_{i,j}$ from $v_i$ to $v_j$, $i<j$, if 
the corresponding lines $\ell_i$ and $\ell_j$ are transverse, 
and an edge $e_{J,i}$ from $v_J$ to $v_i$ if 
$\ell_i \supset F_J$.  See \fullref{fig:nearpencil} 
for an illustration.

Let $m_v$ denote the multiplicity (that is, the degree) of the vertex $v$ 
of $\GA$.  Note that, if $v$ corresponds to the line $L_i$ of 
$\tilde\A$, then $m_v$ is given by the number of lines 
$L_j \in \tilde\A \setminus\{L_i\}$ which intersect $L_i$.  
The graph manifold structure of the boundary manifold 
$M=\partial{N}(\tilde C)$ may be described as follows.  
If $v\in \cV(\GA)$ corresponds to $L_i \in \tilde\A$, 
then the vertex manifold, $M_v$, is given by
\begin{equation} 
\label{eq:vertex manifold}
M_v =\partial U_i \setminus \{\Int(U_j \cap \partial U_i) 
\mid L_j \cap L_i \neq\emptyset \}
\cong S^1 \times \Bigl( \CP^1 \setminus 
\bigcup_{j=1}^{m_v} B_j\Bigr), 
\end{equation}
where $\Int(X)$ denotes the interior of $X$, and the 
$B_j$ are open, disjoint disks.  Note that the boundary 
of $M_v$ is a disjoint union of $m_v$ copies of the torus 
$S^1 \times S^1$.   
The boundary manifold $M$ is obtained by gluing together 
these vertex manifolds along their common boundaries by 
standard longitude-to-meridian orientation-preserving 
attaching maps.

Graph manifolds are often aspherical.  As noted in 
\fullref{ex:boundary pencil}, if $\A$ is a pencil, 
then the boundary manifold of $\A$ is a connected sum 
of $S^1\times S^2$'s, hence fails to be a $K(\pi,1)$--space.  
Pencils are the only line arrangements for which this 
failure occurs.

\begin{prop}[Cohen and Suciu \cite{CS06}]
\label{prop:aspherical} 
Let $\A$ be a line arrangement in $\CP^2$.  The boundary manifold 
$M=M(\A)$ is aspherical if and only if $\A$ is essential, that is, not a pencil.  
\end{prop}

\section{Fundamental group of the boundary}
\label{sec:pi1}

Using the graph manifold structure described in the previous 
section, and a method due to Hirzebruch \cite{Hir}, 
Westlund \cite{We} obtained a presentation for the fundamental 
group of the boundary manifold of a projective line arrangement.  
In this section, we recall this presentation, and use it to obtain 
a minimal presentation.

\subsection{The group of a weighted graph}
\label{subsec:group weighted graph}
Let $\Gamma$ be a loopless graph with $N+1$ vertices.  
Identify the vertex set of $\Gamma$ with $ \set{0,1,\dots,N}$, 
and assume that there is a weight $w_i\in\Z$ given for each 
vertex.  Identify the edge set $\cE$ of $\Gamma$ with a 
subset of $\set{(i,j)\mid 0\le i<j\le N}$ in the obvious manner.  
Direct $\Gamma$ arbitrarily.

We associate a group $G(\Gamma)$ to the weighted graph 
$\Gamma$, as follows.  Let $\cT$ be a maximal tree in $\Gamma$, 
let $\cC=\cE\setminus\cT$, and order the edges in $\cC$.  
Note that $g=\abs{\cC}=b_1(\Gamma)$ is the number of 
(linearly independent) cycles in $\Gamma$.  The group 
$G(\Gamma)$ has presentation
\begin{equation} 
\label{eqn:pres}
G(\Gamma)=
\left\langle
\begin{array}{l}%
 x_0,x_1,\dots,x_N \\[2pt]
\gamma_1,\dots,\gamma_g
\end{array}
\Bigg|
\begin{array}{ll}%
[x_i,x_j^{u_{i,j}}], & (i,j)\in\cE \\ [2pt]
\prod_{j=1}^N x_j^{u_{i,j}}, & 0\le i\le N
\end{array}
\right\rangle,
\end{equation}
where 
\[
u_{i,j}=
\begin{cases}
w_i & \text{if $i=j$,}\\
\gamma_k & \text{if $(i,j)$ is the $k$th element of $\cC$,}\\
\gamma_k^{-1} & \text{if $(j,i)$ is the $k$th element of $\cC$,}\\
1 & \text{if $(i,j)$ or $(j,i)$ belongs to $\cT$,}\\
0 & \text{otherwise}.
\end{cases}
\]
Here $[a,b]=aba^{-1}b^{-1}$, $a^0=1$ is the identity element 
of $G$, and $a^b=b^{-1}ab$ for $b\neq 0$.  Note that if 
$i\neq j$ and $u_{i,j}\neq 0$, then $u_{j,i}=u_{i,j}^{-1}$.

Now let $\A$ be an arrangement of $n+1$ lines in $\CP^2$, 
with associated graph $\GA$, and consider the group $G(\GA)$.  
Recall that the vertices of $\GA$ are in one-to-one 
correspondence with the lines $\set{L_i\mid 0 \le i \le n+r}$ 
of the arrangement $\tilde\A$ in $\widetilde{\CP}{}^2$.  
If $L_i$ is the proper transform of the line $\ell_i\in\A$, 
let $p_i$ denote the number of $0$--dimensional dense edges 
of $\A$ contained in $\ell_i$, and assign the weight 
$w_i=1-p_i$ to the corresponding vertex $v_i$ of $\GA$. 
If $L_i$ is an exceptional line, arising from blowing up 
the dense edge $F_J$ of $\A$, assign the weight $w_J = -1$ 
to the corresponding vertex $v_J$ of $\GA$.  Note that the 
weights of the vertices of $\GA$ are the self-intersection 
numbers of the corresponding lines $L_i$ in $\widetilde{\CP}{}^2$.

\begin{thm}[Westlund \cite{We}] 
\label{thm:pres1}
Let $\A$ be an arrangement of lines in $\CP^2$ with 
boundary manifold $M$. Then the fundamental group of $M$ 
is isomorphic to the group $G(\GA)$ associated to the 
weighted graph $\GA$.
\end{thm}

The presentation provided by this result may be simplified, 
so as to obtain a presentation with $b_1(M)=b_2(M)$ generators 
and relators, realizing $G(\A)=\pi_1(M(\A))$ as a 
commutator-relators group.  The presentation from 
\fullref{thm:pres1} depends on a number of choices: 
the orderings of the lines of $\A$ and the vertices of
$\GA$, the orientation of the edges of $\GA$, and the
choice of maximal tree $\cT$.  As noted by Westlund \cite{We}, 
different choices yield isomorphic groups.  To simplify 
the presentation, we will fix orderings and orientations, 
and work with a specific maximal tree.  Our choice of tree 
will make transparent the relationship between the Betti 
numbers of the boundary manifold $M$ and the complement 
$X$ of $\A$.

\subsection{Simplifying the presentation}
\label{subsec:simpler pres}

Recall that the lines $\{\ell_i\}_{i=0}^n$ of $\A$ are 
ordered.  Designate $\ell_0 \in \A$ as the line at infinity 
in $\CP^2$.  Let $\hat\A$ be the central arrangement in 
$\C^3$ corresponding to $\A\subset \CP^2$, and let $\dA$ 
be the decone of $\hat\A$ with respect to $\ell_0$.  
Incidence with $\ell_0$ gives a partition 
\begin{equation} 
\label{eqn:partition}
\Pi_0=(I_1 \mid I_2 \mid \dots \mid I_f)
\end{equation} 
of the remaining lines of $\A$, where $I_k$ is maximal so that 
$\ell_0 \cap \bigcap_{i\in I_k}\ell_i$ is an edge of $\A$.  
Reorder these remaining lines if necessary to insure
that $I_1=\set{1,\dots,i_1}$, $I_2=\set{i_1+1,\dots,i_2}$, etc., and that 
lines $\ell_i$ transverse 
to $\ell_0$ come last.  In terms of the decone $\dA$ of 
$\A$ with respect to $\ell_0$, this insures that 
members of parallel families of lines in $\dA$ are indexed 
consecutively.

Order the vertices of $\GA$ by $v_{J_1},\dots, v_{J_r}, 
v_1,\dots,v_n,v_0$, where the $v_{J_k}$ are ordered 
lexicographically.  In particular, the vertices 
corresponding to dense edges $F \subset \ell_0$ 
come first.  Recall that the edge $e_{i,j}$ is 
oriented from $v_i$ to $v_j$ if $\ell_i \pitchfork\ell_j$ 
are transverse and $i<j$, and that $e_{J,i}$ is oriented 
from $v_J$ to $v_i$ if the $0$--dimensional dense edge 
$F_J$ is contained in $\ell_i$.

Let $\cT$ be the tree in $\GA$ consisting of the 
following edges:
\[
\cT=
\set{e_{0,i}\mid \ell_0\pitchfork \ell_i}
\cup
\set{e_{J,i}\mid F_J \subset \ell_0\cap \ell_i}
\cup
\set{e_{J,i}\mid F_J\subset \ell_i, \ i=\min{J}}.
\]
It is readily checked that $\cT$ is maximal.  
The edges of $\GA$ not in the tree $\cT$ are 
\[
\cC=
\set{e_{i,j}\mid \ell_i\pitchfork \ell_j,\ 1\le i < j \le n}
\cup
\set{e_{J,i}\mid F_J \subset \ell_i,\ i \neq \min{J},\ 0 \notin J}.
\]
The edges in $\cC$ are in one-to-one correspondence with 
the set $\nbc_2(\dA)$ of pairs of elements of the decone 
$\dA$ which have nonempty intersection and contain 
\emph{no broken circuits}, see Orlik and Terao \cite{OT1}.  It is well 
known that the cardinality of the set $\nbc_2(\dA)$ is
equal to $b_2(X)$, the second Betti number of the 
complement of $\A$.

Now consider the group $G(\A)=G(\GA)$ associated to 
the graph $\GA$.  Denote the generators corresponding 
to the vertices of $\GA$ by 
$x_i$, $0 \le i \le n$, and $x_{J_k}$, $1\le k\le r$, where 
$\set{F_{J_1},\dots,F_{J_r}}$ are the $0$--dimensional dense 
edges of $\A$.  Since the edges of $\cC$ correspond to 
elements $(i,j)\in\nbc_2(\dA)$, we denote the associated 
generators of $G(\A)$ by $\gamma_{i,j}$.  We modify the 
notation of the presentation \eqref{eqn:pres} accordingly, 
writing $R_J$, $u_{J,i}$, $w_{J}$ etc.

\begin{lem}
\label{lem:redundant1}
All commutator relators in the presentation \eqref{eqn:pres} of
$G(\A)=G(\GA)$ involving the generator $x_0$ are redundant.
\end{lem}

\begin{proof}
If $\ell_0 \cap \ell_i$ is a double point of $\A$ for some $i$, 
$1\le i \le n$, then for this $i$, we have the commutator relators 
$[x_p,x_i^{u_{p,i}}]$ for $1\le p<i$ and $\ell_p\pitchfork \ell_i$, 
$[x_i,x_q^{u_{i,q}}]$ for $i<q \le n$ and 
$\ell_i\pitchfork \ell_q$, and 
$[x_{J},x_i^{u_{J,i}}]$ for $F_J \subset \ell_i$.  Here, 
$u_{p,i}=\gamma_{p,i}$, $u_{i,q}=\gamma_{i,q}$, 
$u_{J,i}=1$ if $i=\min{J}$ (by our choice of tree), 
and $u_{J,i}=\gamma_{k,i}$ if $k= \min{J} < i$.
We also have the relator
\[
R_i =
x_{J_1}^{u_{i,J_1}}\ldots x_{J_r}^{u_{i,J_r}}\cdot 
x_1^{u_{i,1}}\ldots x_{i-1}^{u_{i,i-1}}
\cdot x_i^{w_i}\cdot x_{i+1}^{u_{i,i+1}}\ldots 
x_n^{u_{i,n}}\cdot x_0^{u_{i,0}}.
\]
By our choice of tree, we have $u_{i,0}=1$.  
If $\ell_i\cap \ell_j$ is not a double point of $\A$, 
there is no edge joining $v_i$ and $v_j$, and $u_{i,j}=0$.  
Similarly, if $F_J \not\subset \ell_i$, then $u_{i,J}=0$.

Since $u_{p,i}^{}=u_{i,p}^{-1}$ and $u_{J,i}^{}=u_{i,J}^{-1}$, 
the commutator relators $[x_p,x_i^{u_{p,i}}]$ and 
$[x_{J},x_i^{u_{J,i}}]$ are equivalent to $[x_p^{u_{i,p}},x_i]$ and 
$[x_{J}^{u_{i,J}},x_i]$.  It follows that $R_i = a \cdot x_0$, 
where $x_i$ commutes with $a$.  Hence $x_i=x_i^{u_{i,0}}$ commutes with $x_0$.

If $F_J \subset \ell_0$, then $J=\set{i_1,\dots,i_q}$ and $i_1=0$.  
In this instance, we have relators 
$R_{J}= x_{J}^{-1}\cdot x_1^{u_{J,1}}\ldots 
x_n^{u_{J,n}}\cdot x_0^{u_{J,0}}$ and 
$[x_{J},x_{i_p}^{u_{J,i_p}}]$ for $2\le p \le q$.
If $F_J \not\subset \ell_i$, then $u_{J,i}=0$.  
By our choice of tree, $u_{J,i_p}=1$ for 
$1\le p \le q$.  It follows that 
$R_{J}=x_{J}^{-1}\cdot x_{i_2}
\ldots x_{i_q} \cdot x_0$, and $x_{J}$ commutes 
with $x_{i_p}$ for $2\le p\le q$.  Hence $x_{J}=x_J^{u_{J,0}}$ 
commutes with $x_0$.
\end{proof}

Now observe that the relators of type 
$R_{J}=x_{J}^{-1}\cdot \prod_{k=1}^n x_k^{u_{J,k}} \cdot x_0^{u_{J,0}}$ 
may be used to express the generators $x_{J}$ 
in terms of $x_i$, $0\le i \le n$.  If 
$F_J = \ell_{j_1} \cap \dots \cap \ell_{j_q}$
and $j_1=0$, then as noted above, 
$R_{J}=x_{J}^{-1}\cdot x_{j_2}
\ldots x_{j_q} \cdot x_0$.  If $j_1 \ge 1$, then 
$u_{J,k}=0$ for $k\neq j_p$,
$u_{J,j_1}=1$,  
and $u_{J,j_p}=\gamma_{j_1,j_p}$ for $2\le p \le q$.  
So we have
\begin{equation} \label{eqn:big1}
x_{J}=
\begin{cases}
x_{j_2} \ldots x_{j_q}\cdot x_0&\text{if $F_J\subset \ell_0$,}\\
x_{j_1}^{}\cdot x_{j_2}^{\gamma_{j_1,j_2}}\ldots 
x_{j_q}^{\gamma_{j_1,j_q}} &\text{if $F_J\not\subset \ell_0$.}
\end{cases}
\end{equation}
For each $p$, $1\le p\le q$, we have $F_J \subset \ell_{j_p}$ 
and the corresponding commutator relator 
$[x_{J},x_{j_p}^{\gamma_{j_1,j_p}}]$.  In light of
\eqref{eqn:big1}, this may be expressed as
\begin{equation} 
\label{eqn:big2}
\bigl[z_J, \,
x_{j_p}^{\gamma_{j_1,j_p}}\bigr],
\end{equation}
where 
$z_J=x_{j_1}^{}\cdot x_{j_2}^{\gamma_{j_1,j_2}}\ldots 
x_{j_q}^{\gamma_{j_1,j_q}}$ if $j_1 \ge 1$, and 
$z_J=x_{j_2}\ldots x_{j_q} \cdot x_0
=x_0 \cdot x_{j_2}\ldots x_{j_q}$ if $j_1=0$. 

Note that the relator \eqref{eqn:big2} in case $p=1$ (with $\gamma_{j_1,j_1}=1$) is 
a consequence of those for $2\le p\le q$.  Thus, we 
obtain a presentation for $G(\A)$ with generators
$x_i$, $0\le i\le n$, and $\gamma_{i,j}$, $(i,j)\in\nbc_2(\dA)$, the
relators recorded in \eqref{eqn:big2}, together with the relators
$[x_i,x_j^{\gamma_{i,j}}]$, $1\le i<j\le n$, corresponding to double
points $\ell_i\cap \ell_j$ of $\dA$, and $R_i = \prod_{F_J\subset \ell_i}
\smash{x_{J}^{u_{i,J}}}\cdot \prod_{k=1}^n \smash{x_k^{u_{i,k}}} \cdot x_0$, where
$x_{J}$ is given by \eqref{eqn:big1}, the order is irrelevant in the first product, 
and $0\le i \le n$.

\begin{lem} 
\label{lem:redundant2}
If $F_J=\ell_{j_1} \cap \dots \cap \ell_{j_q}$ and $F_J \subset \ell_0$, 
then all the commutator relators recorded in \eqref{eqn:big2} 
are redundant.
\end{lem}

\begin{proof}
We have $j_1=0$ and, by \fullref{lem:redundant1}, the assertion
holds in the case $j_p=0$.  So for $j_p\neq 0$, we must show that the
relator $[x_{J},\,x_{j_p}]$ is a consequence of other relators,
where $x_{J}=x_0\cdot x_{j_2}\ldots x_{j_q}$.

For fixed $j_p{\neq}0$, we have relators
$\bigl[x_i,x_{j_p}^{\gamma_{i,j_p}}\bigr]$ and
$\bigl[x_{j_p},x_k^{\gamma_{j_p,k}}\bigr]$
for $i<j_p<k$, and $\ell_i\pitchfork\ell_{j_p}$, 
$\ell_{j_p}\pitchfork\ell_k$.  The first is equivalent to
$\bigl[\smash{x_i^{\scriptscriptstyle\gamma_{i,j_p}^{-1}}},x_{j_p}\bigr]$.  {F}rom \eqref{eqn:big2}, 
we also have relators $\bigl[x_{j_p},\smash{x_{J_l}^{u_{j_p,J_l}}}\bigr]$ if 
$F_{J_l} \subset \ell_{j_p}$, where $x_{J_l}$ is given 
by \eqref{eqn:big1}, $u_{j_p,J_l}=1$ if $j_p=\min J_l$, and
$u_{j_p,J_l}=\smash{\gamma_{j,j_p}^{-1}}$ if $j_p>j=\min J_l$.  Note that if
$J_l \neq J$, then the word $x_{J_l}$ does not involve the generator
$x_0$.  Additionally, we have the relator $R_{j_p}$, 
which may be expressed as
\[
R_{j_p}=x_{J} \cdot \prod_{J_l\neq J} x_{J_l}^{u_{j_p,J_l}} \cdot 
\prod_{i<j_p}x_i^{\gamma_{i,j_p}^{-1}}\cdot x_{j_p}^{w_{j_p}} \cdot 
\prod_{j_p<k}x_k^{\gamma_{j_p,k}}, 
\] 
where the first product is over all $J_l$ with 
$F_{J_l} \subset \ell_{j_p}$ with $F_{J_l}\not\subset \ell_0$,  
and the last two products are over all $i$, $1\le i<j_p$, 
and $k$, $j_p<k \le n$, for which $\ell_i\pitchfork \ell_{j_p}$ 
and $\ell_{j_p}\pitchfork \ell_k$.

The above commutator relators imply that $R_{j_p}=
x_{J} \cdot a$, where $x_{j_p}$ commutes with $a$. 
Hence $x_{j_p}$ commutes with $x_{J}$.  The result follows.
\end{proof}

\subsection{A commutator-relators presentation}
\label{subsec:comm rels}

There are now $\abs{\nbc_2(\dA)}=b_2(X)$ remaining 
commutator relators: those given by \eqref{eqn:big2} 
corresponding to dense edges  
$F_J = \bigcap_{j \in J} \ell_j$ with $F_J \not\subset \ell_0$, 
and the relators $[x_i,x_j^{\gamma_{i,j}}]$, $1\le i<j \le n$, 
corresponding to double points $\ell_i\pitchfork \ell_j$ of $\dA$.  
Note that all of these commutator relators may be expressed as 
$[z_J,\,x_j^{\gamma_{i,j}}]$, where $\bigcap_{j \in J} \ell_j$ 
is an edge of $\dA$, $i=\min(J)$, and $j \in J\setminus\min(J)$.

There also remain the relators 
\[
R_i = \prod_{F_J\subset \ell_i}
x_{J}^{u_{i,J}} \cdot \prod_{k=1}^n x_k^{u_{i,k}} \cdot x_0^{u_{i,0}},
\]
for $0\le i \le n$.  
We obtain a minimal presentation for $G(\A)$ by eliminating 
the generator $x_0$ using the relator $R_0$.  By our choice 
of tree, this relator is given~by 
\[
R_0=\prod_{F_J \subset \ell_0} x_{J} \cdot x_0^{w_0} \cdot
x_1^{u_{0,1}} \cdots \cdot x_n^{u_{0,n}},
\] 
where $u_{0,i}=1$ if $\ell_0\pitchfork \ell_i$, $u_{0,i}=0$ 
otherwise,  and $x_{J}=x_0\cdot x_{j_2}\ldots x_{j_q}$ if 
$F_J=\ell_0\cap \ell_{j_2}\cap \dots \cap \ell_{j_q}$.  
The chosen ordering of the lines of $\A$ implies that 
$\set{j_2,\dots,j_q}=I_k$, where $(I_1 \mid \dots
\mid I_t)$ is the partition of $\set{1,\dots,n}$ induced by
incidence with $\ell_0$.  Simplifying using the commutation 
relations reveals that 
\begin{equation} 
\label{eq:meridian product}
R_0=x_0 \cdot x_{1} \ldots x_{n}.
\end{equation}
Consequently, we write $x_0=(x_1\ldots x_n)^{-1}$ and 
delete the relation $R_0$.

Now, if $\ell_0\cap \ell_i$ is a double point of $\A$, then
$R_i=Y_i \cdot x_0$, where $Y_i$ is a word in the $x_j$, $j \neq 0$,
and the $\gamma_{i,j}$.  If $\ell_0\cap \ell_i =F_J \in \sD(\A)$, 
then by our ordering of the vertices of $\GA$, 
$R_i = x_{J} \cdot Z_i=x_0 \cdot x_{j_2}\ldots x_{j_q} \cdot Z_i$, 
where $J=\set{0,j_2,\dots,j_q}$.  Conjugating by $x_0$, we 
can write $R_i=Y_i\cdot x_0$, where $Y_i$ is a word as above, 
in this instance as well.  

The next result summarizes the above simplifications.  If
$(i,k) \in \nbc_2(\dA)$, let $F_{I(i,k)}$ be the corresponding 
edge of $\dA$.  For an edge $F_I$ of $\dA$, with $i=\min{I}$, 
and $j \in I \setminus \min{I}$, let $\gamma_{I,j}=\gamma_{i,j}$.
If $\ell_0 \cap \ell_p \cap \dots \cap \ell_q$ is an edge of 
$\A$, set $\zeta_{0,j}= x_p \ldots x_q$ for each $j$, 
$p \le j \le q$.  Note that if $\ell_0$ and $\ell_j$ are
transverse, then $\zeta_{0,j}=x_j$.

\begin{prop} 
\label{prop:THEpres}
The fundamental group of the boundary manifold $M$ of $\A$ has presentation
\begin{equation*}
G(\A)=
\left\langle
\begin{array}{ll}%
x_j,& 1\le j\le n\\[2pt]
\gamma_{i,k},& (i,k)\in \nbc_2(\dA)
\end{array}
\Bigg|
\begin{array}{ll}%
R_j,& 1\le j\le n\\[2pt]
R_{i,k},& (i,k)\in \nbc_2(\dA)
\end{array}
\right\rangle,
\end{equation*}
where
\[
R_j = \zeta_{0,j} \cdot \hskip -6pt
\prod_{\atop{F_I\in \sD(\dA)}{j\in I \setminus \min{I}}} 
(\gamma_{I,j}^{} z_I^{}\gamma_{I,j}^{-1}x_j^{-1}) \cdot
\hskip -6pt
\prod_{\atop{F_I\in \sD(\dA)}{j = \min{I}}} 
(x_j^{-1}z_I^{}) \cdot
 \hskip -6pt
\prod_{\atop{\ell_i \pitchfork \ell_j}{1\le i<j}} 
x_i^{\gamma_{i,j}^{-1}} \cdot
 \hskip -6pt
\prod_{\atop{\ell_j \pitchfork \ell_k}{j<k\le n}} 
x_k^{\gamma_{j,k}^{}} \cdot
(x_1\ldots x_n)^{-1}
\]
and 
\[
R_{i,k} = [z_{I(i,k)}^{},\,x_k^{\gamma_{i,k}}].
\]
\end{prop}

\begin{proof}
It follows from the preceding discussion that the group 
$G(\A)$ has such a presentation with the relators $R_{i,k}$ 
as asserted.  So it is enough to show that the relators $R_j$ 
admit the above description.

Fix $j$, $1\le j \le n$, and consider the line $\ell_j$ of $\A$.  
Assume that
\begin{equation} 
\label{eqn:jdata}
\begin{aligned}
\hfill\mathrm{(i)}\quad & \text{$j \in J$, where $J=[p,q]$ and 
$\ell_0 \cap \ell_p \cap \dots \cap \ell_q$ is an edge of $\A$;} \\
\hfill\mathrm{(ii)}\quad & \text{$j \in J_t \setminus \min{J_t}$ 
for $1\le t \le a$ and 
$F_{J_t}$ is a dense edge of  $\dA$;} \\
\mathrm{(iii)}\quad & \text{$j =\min{K_t}$ for $1\le t \le b$ and 
$F_{K_t}$ is a dense edge of  $\dA$;} \\
\hfill\mathrm{(iv)}\quad & \text{$\ell_j \pitchfork \ell_{i_t}$ for 
$1\le t \le c$ and $1\le i_t <j$;
and } \\
\hfill\mathrm{(v)}\quad & \text{$\ell_j \pitchfork \ell_{k_t}$ 
for $1\le t \le d$ and $j < k_t \le n$.}
\end{aligned}
\end{equation}
Note that $\ell_j$ contains either $a+b$ or $a+b+1$ dense edges of $\A$, 
depending on whether $\ell_j$ is transverse to $\ell_0$ or
not.  Consequently, the weight of the vertex $v_j\in\GA$ is
\[
w_j=\begin{cases}
1-a-b& \text{if $\ell_j \pitchfork \ell_0$,}\\
-a-b & \text{otherwise.}
\end{cases}
\]
With these data, the preceding discussion and our conventions
regarding the graph $\GA$ and the group $G(\A)$ imply that the 
relator $R_j$ is given by
\[
R_j = \zeta_{0,j} \cdot 
\prod_{t=1}^a z_{J_t}^{\gamma_{j,J_t}}
\cdot
\prod_{t=1}^b z_{K_t}^{\gamma_{j,K_t}}
\cdot
\prod_{t=1}^c x_{i_t}^{\gamma_{j,i_t}}
\cdot
x_j^{-a-b} 
\cdot 
\prod_{t=1}^d x_{k_t}^{\gamma_{j,k_t}}
\cdot
x_0.
\]
The commutator relators $R_{i,k}$ imply that $x_j$ commutes with each
of $z_{J_t}^{\gamma_{j,J_t}}$, $z_{K_t}^{\gamma_{j,K_t}}$,
$x_{i_t}^{\gamma_{j,i_t}}$, $x_{k_t}^{\gamma_{j,k_t}}$ for all
relevant $t$.  Furthermore, $\gamma_{j,J}^{}=\gamma_{J,j}^{-1}$ if 
$j \in J \setminus \min{J}$, $\gamma_{j,K}=1$ if $j=\min{K}$, and
$\gamma_{j,i}^{}=\gamma_{i,j}^{-1}$ if $i<j$.  Using these facts, the
relator $R_j$ may be expressed as
\[
R_j = \zeta_{0,j} \cdot 
\prod_{t=1}^a (z_{J_t}^{\gamma_{J_t,j}^{-1}}\cdot x_j^{-1})
\cdot
\prod_{t=1}^b (x_j^{-1}\cdot z_{K_t}^{})
\cdot
\prod_{t=1}^c x_{i_t}^{\gamma_{i_t,j}^{-1}}
\cdot
\prod_{t=1}^d x_{k_t}^{\gamma_{j,k_t}}
\cdot
x_0.
\]
Recalling that $x_0=(x_1\ldots x_n)^{-1}$, this is easily seen to be
equivalent to the expression given in the statement of the
Proposition.
\end{proof}

\begin{remark} 
\label{rem:randell} 
If $(i,k)\in\nbc_2(\dA)$ and $I=I(i,k)=\set{i_1,\dots,i_q}$,
the relators $[z_I,\,x_{i_p}^{\gamma_{i_1,i_p}}]$, $2\le p \le q$, are
equivalent to the family
$[x_{i_1}^{},\,x_{i_2}^{\gamma_{i_1,i_2}},\dots,x_{i_q}^{\gamma_{i_1,i_q}}]$ 
of ``Randell relations'' familiar from presentations of the fundamental 
group of the complement of an arrangement.
\end{remark}

\begin{cor} 
\label{cor:commrels}
The group $G(\A)$ is a commutator-relators group.
\end{cor}
\begin{proof}
By \fullref{prop:THEpres}, the group $G(\A)=\pi_1(M)$ 
admits a presentation with $b_1=b_1(M)$ generators.  The 
conclusion follows from this, together with the fact that $H_1(M)$ 
is free abelian of rank $b_1$, see Matei and Suciu \cite[Proposition 2.7]{MS00}.
\end{proof}

\begin{remark} 
This result may also be established directly, by showing that each relator 
$R_j$ is a product of commutators.  Using the Randell relations 
noted above, one can show that 
$R_j=x_{\rho_{i,1}}^{v_{i,1}}\ldots x_{\rho_{i,n}}^{v_{i,n}} \cdot
x_0$, where $\set{\rho_{i,1},\dots,\rho_{i,n}}$ is a permutation of
$[n]$ and $v_{p,q}$ is a word in the generators $\gamma_{i,j}$.   
This may be expressed as a  product of commutators using 
the fact that $x_0=(x_{1} \ldots x_{n})^{-1}$.
\end{remark}

\subsection{Some computations}  
\label{subsec:gp and np}

We conclude this section with a few examples illustrating 
how the presentation from \fullref{prop:THEpres} 
works in practice. 

\begin{example} 
\label{example:near-pencil pres}
Let $\A$ be a near-pencil of $n{+}1\ge4$ lines,
with defining polynomial $Q=x_0(x_1^n-x_2^n)$ 
and boundary manifold $M$.
The graph $\GA$ has
vertices $v_0,v_1,\dots,v_n$ corresponding to the lines, and
one more vertex $v_{n+1}=v_F$ corresponding to the multiple
point $F=\ell_1\cap \dots\cap \ell_n$.  The
weights of the vertices are $w_0=1$, $w_1=\dots=w_n=0$, and
$w_{n+1}=-1$.  The  edge set is
$\cE$ consists of edges $e_{0,i}$ and $e_{i,n+1}$ for $1\le i \le n$.  
Fix the maximal tree
$\cT=\{e_{0,1},\dots,e_{0,n},e_{1,n+1}\}$, 
indicated by dashed edges in \fullref{fig:nearpencil}.

By \fullref{prop:THEpres}, the
fundamental group of $M$ 
has presentation
\[
G(\A)=\langle
x_1,\ x_j,\ \gamma_{1,j} \mid
z^{} \zeta^{-1},\ x_j^{} \gamma_{1,j}^{} z_{}^{} \gamma_{1,j}^{-1}
x_j^{-1} \zeta_{}^{-1},\ [z,x_j^{\gamma_{1,j}}]
\rangle,
\]
where $z=x_1^{}\cdot x_2^{\gamma_{1,2}} \ldots 
x_n^{\gamma_{1,n}}$, $\zeta=x_1\cdot x_2\ldots x_n$, and $2 \le j \le n$.

The elements $\zeta,x_2,\dots,x_n,\gamma_{1,2},\dots,\gamma_{1,n}$
generate the group $G(\A)$, and it is readily checked that $\zeta$ is
central.  Also, conjugating the relator $R_1$ by $x_1$ yields
\[
[\gamma_{1,2}^{-1},x_2^{}] \cdot
x_2^{}[\gamma_{1,3}^{-1},x_3^{}]x_2^{-1}\ldots
(x_2^{}\ldots x_{n-1}^{})[\gamma_{1,n}^{-1},x_{n}^{}](x_2^{}\ldots
x_{n-1}^{})^{-1}.
\]
It follows that $G(\A)$ is isomorphic to the direct product
of a cyclic group  $\Z=\langle c\rangle$ with a genus $n-1$
surface group 
\[
\pi_1(\Sigma_{n-1})=
\langle g_1,\dots,g_{2n-2}\mid [g_1,g_2]\ldots
[g_{2n-3},g_{2n-2}]\rangle.
\]
An explicit isomorphism $\Z \times
\pi_1(\Sigma_{n-1})\xrightarrow{\simeq} G(\A)$ is given by
\[
c\mapsto \zeta,\quad
g_i \mapsto
\begin{cases}
x_2^{}\ldots x_{k}^{}\cdot \gamma_{1,k+1}^{-1} \cdot (x_2^{}\ldots
x_{k}^{})^{-1},&\text{if $i=2k-1$,}\\
x_2^{}\ldots x_{k}^{}\cdot x_{k+1}^{} \cdot (x_2^{}\ldots
x_{k}^{})^{-1},&\text{if $i=2k$}.
\end{cases}
\]
\end{example}

\begin{example} 
\label{ex:general position}
Let $\A$ be an arrangement of $n+1$ lines in general position.
The graph $\GA$ is the complete graph on $n+1$ vertices.  
Here, there are no $0$--dimensional dense edges and all vertices 
have weight $1$.  
\begin{figure}%
\subfigure{%
\label{fig:gen}%
\begin{minipage}[t]{0.4\textwidth}
\setlength{\unitlength}{16pt}
\begin{picture}(4.5,4.8)(-3.3,-0.8)
\put(0,0){\line(1,1){4}}
\put(-0.9,1.4){\line(1,0){5.8}}
\put(0,4){\line(1,-1){4}}
\put(0.5,-0,5){\line(0,1){5}}
\put(1.2,-0,5){\makebox(0,0){$\ell_0$}}
\put(4.4,-0.5){\makebox(0,0){$\ell_1$}}
\put(5.5,1.4){\makebox(0,0){$\ell_2$}}
\put(4.4,4.5){\makebox(0,0){$\ell_3$}}
\end{picture}
\end{minipage}
}
\subfigure{%
\label{fig:gengraph}%
\begin{minipage}[t]{0.4\textwidth}
\setlength{\unitlength}{18pt}
\begin{picture}(5,4.8)(-2.3,-2.45)
\xygraph{!{0;<15mm,0mm>:<0mm,13mm>::}
[]*D(2.5){v_1}*{\,\disc}
(
-|{\gamma_{1,2}}[ddr]*L(1.8){v_2}*{\disca}(-|{\gamma_{2,3}}[ll]),
-|{\gamma_{1,3}}[ddl]*R(2){v_3}*{\disc},
-@{--}[d]*U(3.3){v_0}*{\,\disc}
(-@{--}[dr],-@{--}[dl])
)
}
\end{picture}
\end{minipage}
}
\caption{A general position arrangement and 
its associated graph}
\label{fig:generic}
\end{figure}
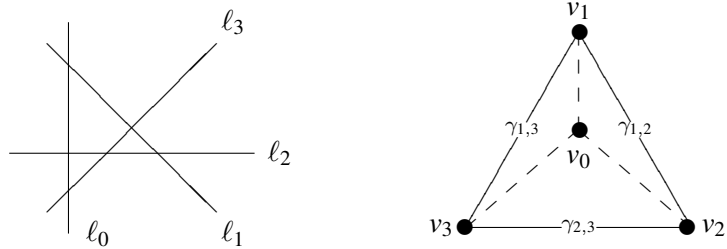

Using the maximal tree $\cT=\set{e_{0,i}\mid 1\le i\le n}$ 
(indicated by dashed edges in \fullref{fig:generic}),
\fullref{prop:THEpres} yields a presentation for $G(\A)$
with generators $x_i$ ($1\le i \le n$) and $\gamma_{i,j}$
($1\le i<j \le n$), and relators 
\begin{align*}
R_j&=x_j \cdot 
x_1^{\gamma_{1,j}^{-1}}\ldots x_{j-1}^{\gamma_{j-1,j}^{-1}}
\cdot
x_{j+1}^{\gamma_{j,j+1}^{}} \ldots x_{n}^{\gamma_{j,n}^{}}
\cdot x_n^{-1} \ldots x_1^{-1}
&& (1\le j \le n),
\\
R_{i,j}&=[x_i^{},x_j^{\gamma_{i,j}}] 
&& (1\le i< j \le n).
\end{align*}
\end{example}

\section{Twisted Alexander polynomials}
\label{sec:AlexPolys}

A finitely generated module $K$ over a Noetherian ring $R$ 
admits a finite presentation
$$R^r \xrightarrow{~~\psi~~}  R^s \longrightarrow K \longrightarrow 0.$$ 
Let $E_i(K)$ denote the $i$th elementary ideal of $K$, the ideal of 
$R$ generated by the codimension $i$ minors of the matrix 
$\psi$.  It is well known that the elementary ideals do not depend on 
the choice of presentation, so are invariants of the module $K$.  

Let $\Lambda=\bF[t_1^{\pm 1},\dots, t_n^{\pm 1}]$ be the  
ring of Laurent polynomials in $n$ variables over a field $\bF$. 
Since $\Lambda$ is a unique factorization domain, there is a 
unique minimal principal ideal that contains the elementary ideal 
$E_0(K)$.  Define the order, $\ord(K)$, of the module $K$ to be 
a generator of this principal ideal.  Note that $\ord(K)$ is defined 
up to multiplication by a unit in $\Lambda$, which necessarily 
is of the form $c t_1^{l_1} \ldots t_n^{l_n}$, for some $l_i\in \Z$ 
and $c \in \bF^*$. 

Now let $G$ be a group of type $FL$, and 
$\alpha\colon G \to H$ a homomorphism to a finitely 
generated, free abelian group.  Note that if $\rank(H)=n$, 
then $\bF[H] \cong \Lambda$.  Let $\phi\colon G \to \GL_k(\bF)$ 
be a representation.  With these data, the vector space 
$\Lambda^k_{\phi,\alpha} = \bF^k \otimes_\bF \Lambda$ 
admits the structure of a (left) $G$--module: if $\gamma \in G$ 
and $v \otimes q \in \smash{\Lambda^k_{\phi,\alpha}}$, then
\[
\gamma \cdot (v \otimes q) = (\phi(\gamma) v) \otimes (\alpha(\gamma)q).
\]
Following Kirk and Livingston \cite{KL}, define the twisted Alexander modules 
of $G$ (with respect 
to $\alpha$ and $\phi$) to be the homology groups 
of $G$ with coefficients in $\Lambda^k_{\phi,\alpha}$: 
if $C_*(G)$ is a finite, free resolution of $\Z$ over $\Z{G}$, then
\begin{equation}
\label{eq:twistmodule}
H_i(G;\Lambda^k_{\phi,\alpha})=H_i(C_*(G) \otimes_{\Z{G}} 
\Lambda^k_{\phi,\alpha}).
\end{equation} 
Note that $H_i(G;\Lambda^k_{\phi,\alpha})$ carries the structure 
of a (finitely generated) right $\Lambda$--module. 
Define the twisted Alexander polynomial $\Delta^{\phi,\alpha}_i(G)$ 
to be the order of this module:
\begin{equation}
\label{eq:twistpoly}
\Delta^{\phi,\alpha}_i(G)=\ord\bigl(H_i(G;\Lambda^k_{\phi,\alpha})\bigr).
\end{equation}
If $\theta\colon G\surj G'$ is an epimorphism, $\alpha=\alpha'\circ \theta$, 
and $\phi=\phi'\circ \theta$, then $\Delta_1^{\phi',\alpha'}(G')$ divides 
$\Delta_1^{\phi,\alpha}(G)$, see Kitano, Suzuki and Wada \cite{KSW}. 

In the case where $\alpha\colon G \surj H_1(G)/\Tors(H_1(G))$ 
is the projection onto the maximal torsion-free abelian quotient, 
we suppress $\alpha$ and write simply $\Lambda^k_\phi$ and 
$\Delta_i^\phi(G)$.  Note that if $\phi\colon G \to \GL_1(\bF)$ is 
the trivial representation, then $\Delta_1^\phi(G)$ is the classical 
Alexander polynomial $\Delta(G)$.  Up to a monomial change of 
variables, $t_i \mapsto t_1^{a_{i,1}}\ldots t_n^{a_{i,n}}$, 
where $(a_{i,j}) \in \GL_n(\Z)$, this Laurent polynomial is an 
invariant of the isomorphism type of the group $G$. 
In what follows, we will focus our attention on the case $\bF=\C$.

\begin{lem} \label{lem:free abelian}
Let $G$ be a finitely generated free abelian group, and 
$\phi\colon G \to \GL_k(\C)$ a representation.  Then the twisted 
Alexander module $H_i(G;\Lambda^k_\phi)$ vanishes for $i \ge 1$, 
and $\ord \bigl( H_0(G;\Lambda^k_\phi) \bigr)= 1$.
\end{lem}

\begin{proof}
Let $n=\rank(G)$. Denote the generators of $G$ 
by $t_1,\dots,t_n$, and identify $\C [G] \cong \Lambda$.   

The proof is by induction on $k$.  If $k=1$, the 
chain complex $C_*(G)\otimes_{\Z{G}} \Lambda^1_\phi$ may 
be realized as the standard Koszul complex in the variables 
$z_i = \phi(t_i)\cdot t_i$. Consequently, $H_i(G;\Lambda^1_\phi) = 
H_i(C_*(G)\otimes_{\Z{G}}\Lambda^1_\phi) = 0$ for $i \ge 1$, and 
$H_0(G;\Lambda^1_\phi) =\C$ has order $1$.

Suppose $k>1$.  Since $G$ is 
abelian, the automorphisms $\phi(t_i) \in \GL_k(\C)$, 
$1\le i \le n$, all commute.  Consequently, they have 
a common eigenvector, say $v$.  Let $\lambda_i$ be the 
eigenvalue of $\phi(t_i)$ with eigenvector $v$, and let 
$\{w_1,\dots,w_{k-1}\}$ be a basis for $\langle v \rangle^\perp$.  
With respect to the basis $\{v,w_1,\dots,w_{k-1}\}$ for $\C^k$, 
the matrix $A_i$ of $\phi(t_i)$ is of the form
\[
A_i = \begin{pmatrix} \lambda_i & * \\ 0 & \bar{A}_i \end{pmatrix}
\]
where $\bar{A}_i$ is an invertible $(k-1)\times (k-1)$ matrix.  
Define representations $\phi'\colon G \to \C^*$ and 
$\phi''\colon G \to \GL_{k-1}(\C)$ by $\phi'(t_i) = 
\lambda_i$ and $\phi''(t_i) = \bar{A}_i$.
Then we have a short exact sequence of $G$--modules
\[
\xymatrix{
0 \ar[r] & \Lambda^1_{\phi'} \ar[r] & \Lambda^k_\phi \ar[r] 
& \Lambda^{k-1}_{\phi''} \ar[r] & 0
},
\]
and a corresponding long exact sequence in homology
\[
\xymatrix{
\dots\ar[r] & H_i(G;\Lambda^1_{\phi'}) \ar[r] &H_i(G;\Lambda^k_{\phi}) 
\ar[r] & H_i(G;\Lambda^{k-1}_{\phi''})\ar[r] & \dots
}.
\]
Using this sequence, the case $k=1$, and the inductive hypothesis, 
we conclude that $H_i(G;\Lambda^k_\phi)=0$ for $i \ge 1$, and 
that $\ord H_0(G;\Lambda^k_\phi) = 1$.
\end{proof}

Let $\Gamma$ be a connected, directed graph, 
and let  $\cV=\cV(\Gamma)$ and $\cE=\cE(\Gamma)$ 
denote the vertex and edge sets of $\Gamma$.  A graph of groups is 
such a graph, together with vertex groups $\{G_v \mid v \in \cV\}$, 
edge groups $\{G_e \mid e \in \cE\}$, and monomorphisms 
$\theta_0\colon G_e \to G_v$ and $\theta_1\colon G_e \to G_w$ 
for each directed edge $e=(v,w)$.  Choose a maximal tree $T$ 
for $\Gamma$.  The fundamental group $G=G(\Gamma)$ 
(relative to $T$) is the group generated by the vertex groups 
$G_v$ and the edges $e$ of $\Gamma$ not in $T$, with the 
additional relations
$e\cdot \theta_1(x) =  \theta_0(x) \cdot e$, for $x\in G_e$ if 
$e \in \Gamma\setminus T$, and 
$\theta_1(y) =  \theta_0(y)$, for $y\in G_e$ if $e \in T$.

\begin{thm} 
\label{thm:graph of groups}
Let $(\Gamma, \{G_e\}_{e\in\cE(\Gamma)}, 
\{G_v\}_{v\in\cV(\Gamma)})$ be a graph of groups, 
with fundamental group $G$, vertex groups of type 
$FL$, and free abelian edge groups.  Assume that the inclusions 
$G_e \hookrightarrow G$ induce monomorphisms in homology.  
If $\phi \colon G \to \GL_k(\C)$ is a representation, then 
\begin{romenum}
\item  \label{gg1}
$H_i(G;\Lambda^k_\phi) = \bigoplus_{v\in \cV} 
H_i(G_v;\Lambda^k_\phi)$ for $i \ge 2$, and
\item \label{gg2}
$\ord\bigl( H_1(G;\Lambda^k_\phi) \bigr) = 
\ord\bigl( \bigoplus_{v\in \cV} H_1(G_v;\Lambda^k_\phi)\bigr)$.
\end{romenum}
\end{thm}

\begin{proof}
For simplicity, we will suppress the coefficient module  
$\Lambda^k_\phi$ for the duration of the proof. 
Given a graph of groups, there is a Mayer--Vietoris sequence
\begin{equation*}
\xymatrixcolsep{14pt}
\xymatrix{
\dots\! \ar[r]& \bigoplus_{e \in \cE} H_i(G_e) 
\ar[r]& \bigoplus_{v \in \cV}  H_i(G_v)
\ar[r]& H_i(G) 
\ar[r]^(.34){\partial}& \bigoplus_{e \in \cE} 
H_{i-1}(G_e)
\ar[r]& \!\dots
}
\end{equation*}
see Brown \cite[Section VII.9]{brown}.  Since the edges groups are free 
abelian and the inclusions $G_e \hookrightarrow G$ induce 
monomorphisms in homology, we may apply 
\fullref{lem:free abelian} to conclude that 
$H_i(G_e)=0$ for all $i \ge 1$.  
Assertion \eqref{gg1} follows. 

\fullref{lem:free abelian}  also implies that
$\ord\bigl(H_0(G_e)\bigr)=1$, 
for each $e \in \cE$.  Consequently, 
$\ord\bigl(\bigoplus_{e\in\cE}H_0(G_e)\bigr)=1$.  
The above Mayer--Vietoris sequence reduces to
\[
\xymatrix{
0 \ar[r] &\bigoplus_{v \in \cV}  H_1(G_v)
\ar[r] & H_1(G) \ar[r]^(.38){\partial} & 
\bigoplus_{e \in \cE} H_{0}(G_e)  \ar[r] & \dots
}.
\]
{F}rom this, we obtain a short exact sequence
\[
\xymatrix{
0 \ar[r] & \bigoplus_{v \in \cV}  H_1(G_v)
\ar[r] & H_1(G) \ar[r]^(.5){\partial} &
\Image(\partial) \ar[r] & 0
}.
\]
Since $\Image(\partial)$ is a submodule of $\bigoplus_{e \in \cE} 
H_{0}(G_e)$, and the latter has order $1$, we have 
$\ord\bigl(\Image(\partial)\bigr)=1$ as well.   
Assertion \eqref{gg2} follows. 
\end{proof}

\section{Alexander polynomials of line arrangements} 
\label{sec:alex poly arr}
Let $\A=\{\ell_0,\dots ,\ell_n\}$ be an arrangement of 
$n+1$ lines in $\CP^2$, with boundary manifold $M$.  
Since $M$ is a graph manifold, the fundamental group 
$G=\pi_1(M)$ is the fundamental group of a graph of 
groups. Recall from \fullref{sec:bdry} that, in the 
graph manifold structure, the vertex manifolds are of 
the form $M_v \cong S^1 \times \bigl( \CP^1 \setminus 
\bigcup_{j=1}^{m} B_j\bigr)$, where the $B_j$ are disjoint 
disks and $m$ is the multiplicity (degree) of the vertex 
$v$ of $\GA$, and these vertex manifolds are glued 
together along tori.  Consequently, the vertex groups 
are of the form $\Z \times F_{m-1}$, and the edge groups 
are free abelian of rank $2$.  

The edge groups are generated by meridian loops 
about the lines $L_i$ of $\tilde\A$ in $\widetilde{\CP}{}^2$.  
In terms of the generators $x_i$ of $G$, these generators 
are of the form $x_i^y$ or $x_{i_1}^{y_1}\ldots x_{i_k}^{y_k}$ 
if $L_i$ is the proper transform of $\ell_i\in\A$ or $L_i$ is 
the exceptional line arising from blowing up the dense edge 
$F_I$ of $\A$, where $I=\set{i_1,\dots,i_k}$.
By \eqref{eq:meridian product}, $x_0 x_1 \ldots x_n=1$ 
in $G$.  This fact may be used to check that the inclusions 
of the edge groups in $G$ induce monomorphisms in homology.  
Therefore, \fullref{thm:graph of groups} may be applied 
to calculate twisted Alexander polynomials of $G$.  We first 
record a number of preliminary facts.

\begin{lem} 
\label{lem:hopf alex}
Let $G = \Z \times F_{m-1}$, and let $\phi\colon G \to \GL_k(\C)$ 
be a representation.  Then the twisted Alexander polynomial 
$\Delta^\phi_1(G)$ is given by
\[
\Delta^\phi_1(G) = \bigl[ p(A,t) \bigr]^{m-2},
\]
where $t$ is the image of a generator $z$ of the center $\Z$ 
of $G$ under the abelianization map, and $p(A,t)$ is the 
characteristic polynomial of the automorphism $A=\phi(z)$ 
in the variable $t$.  In particular, the classical Alexander 
polynomial is $\Delta(G)=(t-1)^{m-2}$.
\end{lem}
\begin{proof}
Write $G=\Z\times F_{m-1} = \langle z,y_1,\dots,y_{m-1} \mid 
[z,y_1], \dots, [z,y_{m-1}]\rangle$.  Applying the Fox calculus 
to this presentation yields a free $\Z{G}$--resolution of $\Z$, 
\[
\xymatrix{
(\Z{G})^{m-1} \ar^{\partial_2}[r]& (\Z{G})^m \ar^(.55){\partial_1}[r] 
& \Z{G} \ar^{\epsilon}[r] & \Z \ar[r] & 0
},
\]
where $\epsilon\colon \Z{G} \to \Z$ is the augmentation map, and 
the matrices of $\partial_1$ and $\partial_2$ are given by 
$[\partial_1]=\begin{pmatrix} z-1 & y_1-1 & \cdots 
& y_{m-1}-1\end{pmatrix}^\top$ and
\[
[\partial_2] = \begin{pmatrix}
1-y_1 & z-1 & 0 & \cdots & 0 \\
1-y_2 & 0 & z-1 & \cdots & 0 \\
\vdots & \vdots & \vdots & \ddots & \vdots \\
1-y_{m-1} & 0 & 0 & \cdots & z-1
\end{pmatrix}
\]
A calculation with this resolution yields the result.
\end{proof}

Let $\GA$ denote the graph underlying the graph manifold 
structure of the boundary manifold $M$ of the line arrangement 
$\A=\{\ell_i\}_{i=0}^n$ in $\CP^2$ and the graph 
of groups structure of the fundamental group $G=\pi_1(M)$.  
For a vertex $v$ of $\GA$ with multiplicity $m_v$, in the 
identification $G_v = \Z \times F_{m_v-1}$ of the vertex 
groups of $G$, the center $\Z$ of $G_v$ is generated by 
$z_v$, an meridian loop about the corresponding line 
$L_i$ of $\tilde\A$.  Denoting the images of the 
generators $x_i$ of $G$ under the abelianization 
$\alpha\colon G \to G/G'$  by $t_i$, there is a 
choice of generator $z_v$ so that 
\[
\alpha(z_v)=t_v=\begin{cases}
t_i & \text{if $v=v_i$, $0 \le i \le n$,}\\
t_{i_1} \ldots t_{i_k} & \text{if $v=v_I$, where 
$I=\set{i_1,\dots,i_k}$ and $F_I\in\sD(\A)$.}
\end{cases}
\]
If $I=\set{i_1,\dots,i_k}$, we subsequently write 
$t_I= t_{i_1} \ldots t_{i_k}$.

\fullref{thm:graph of groups} and \fullref{lem:hopf alex} 
yield the following result.

\begin{thm} 
\label{thm:alex poly arr}
Let $\A$ be an essential line arrangement in $\CP^2$, let $\GA$ 
be the associated graph, and let $G$ be the fundamental 
group of the boundary manifold $M$ of $\A$.  If 
$\phi\colon G \to \GL_k(\C)$ is a representation, 
then the twisted Alexander polynomial $\Delta^\phi_1(G)$ 
is given~by
\[
\Delta^\phi_1(G) = \prod_{v \in \cV(\GA)} \bigl[ p(A_v,t_v) \bigr]^{m_v-2}, 
\]
where $t_v$ is the image of a generator of the center $\Z$ 
of $G_v$ under the abelianization map, and $p(A_v,t_v)$ is 
the characteristic polynomial of the automorphism $A_v=\phi(z_v)$ 
in the variable $t_v$.  In particular, the classical Alexander 
polynomial of $G$ is 
\[
\Delta(G) = \prod_{v \in \cV(\GA)} (t_v-1)^{m_v-2}.
\]
\end{thm}

\begin{remark} 
By gluing formulas of Meng and Taubes \cite{MT} and Turaev \cite{Tu}, 
with appropriate identifications, Milnor torsion is multiplicative 
when gluing along tori.  Since Milnor torsion coincides with the 
Alexander polynomial for a $3$--manifold $M$ with $b_1(M)>1$, the
calculation of $\Delta(G)$ in \fullref{thm:alex poly arr} 
above may alternatively be obtained using these gluing formulas, 
see Vidussi \cite[Lemma~7.4]{Vi}.

The above formula for $\Delta(G)$ is also reminiscent of the 
Eisenbud--Neumann formula for the Alexander polynomial $\Delta_L(t)$ 
of a graph (multi)-link $L$, see Eisenbud and Neumann
\cite[Theorem~12.1]{EN}.  For example, if $L$ is the $n$--component
Hopf link (that is, the singularity link of a pencil of $n\ge 2$ lines),
then $\Delta_L(t)=(t_1\ldots t_n-1)^{n-2}$.
\end{remark}

Recall from \eqref{eq:meridian product} that the meridian generators 
$x_i$ of $G$ corresponding to the lines $\ell_i$ of $\A$, $0 \le i \le n$, 
satisfy the relation $x_0 x_1 \ldots x_n=1$.  Consequently, 
$t_0 t_1 \ldots t_n=1$ and the twisted Alexander polynomial 
$\Delta_1^\phi(G)$ may be viewed as an element of 
$\Lambda=\C[t_1^{\pm 1},\dots,t_n^{\pm 1}]$.  In particular, in the 
classical Alexander polynomial, if $I = \{0\} \cup J$, then 
$t_I - 1 \doteq t_{[n]\setminus J} - 1$, since Alexander polynomials 
are defined up to multiplication by units.  In what follows, we make 
substitutions such as these without comment.

In light of \fullref{thm:alex poly arr}, we focus on the classical 
Alexander polynomial for the remainder of this section.

\begin{example} 
\label{ex:falk}
In \cite{Fa}, Falk considered a pair of arrangements whose 
complements are homotopy equivalent, even though the two 
intersection lattices are not isomorphic. In this example, 
we analyze the respective boundary manifolds.   

The Falk arrangements $\cF_1$ and $\cF_2$ have defining polynomials
\begin{align*}
Q(\cF_1)&=x_0(x_1+x_0)(x_1-x_0)(x_1+x_2)x_2(x_1-x_2) \\
\text{and}\qquad
Q(\cF_2)&=x_0(x_1+x_0)(x_1-x_0)(x_2+x_0)(x_2-x_0)(x_2+x_1-x_0).
\end{align*}
These arrangements, and the associated graphs, are depicted in 
Figures \ref{fig:falk1} and~\ref{fig:falk2}.

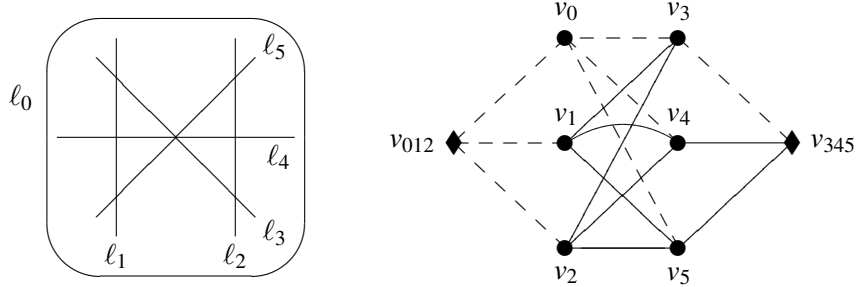
\begin{figure}%
\subfigure{%
\label{fig:f1}%
\begin{minipage}[t]{0.4\textwidth}
\setlength{\unitlength}{15pt}
\begin{picture}(5,5.8)(-3.6,-1.5)
\put(2,2){\oval(6.5,6)[t]}
\put(2,2){\oval(6.5,7)[b]}
\put(0,0){\line(1,1){4}}
\put(-1,2){\line(1,0){6}}
\put(0,4){\line(1,-1){4}}
\put(0.5,-0,5){\line(0,1){5}}
\put(3.5,-0,5){\line(0,1){5}}
\put(-1.9,3){\makebox(0,0){$\ell_0$}}
\put(0.5,-1){\makebox(0,0){$\ell_1$}}
\put(3.5,-1){\makebox(0,0){$\ell_2$}}
\put(4.5,-0.4){\makebox(0,0){$\ell_3$}}
\put(4.6,1.5){\makebox(0,0){$\ell_4$}}
\put(4.5,4.3){\makebox(0,0){$\ell_5$}}
\end{picture}
\end{minipage}
}
\subfigure{%
\label{fig:f1graph}%
\begin{minipage}[t]{0.41\textwidth}
\setlength{\unitlength}{18pt}
\begin{picture}(5,5.3)(-0.5,-2.8)
\xygraph{!{0;<15mm,0mm>:<0mm,14mm>::}
!~:{@{-}|@{~}}
[]*D(3){v_0}*{\disc}
(
-@{--}[dr]
,-@{--}[dl]*R(1.8){v_{012}}*-{\blacklozenge}
(-@{--}[dr]*U(3){v_2}*{\disc}(-[r],-[uur],-[ur])
,-@{--}[r])
,[d]*D(3){v_1}*{\disc}(-[ur],-[dr])
,[dr]*D(3){v_4}*-{\disc}
(
-@/^-7pt/[l]
,[r]*L(1.8){v_{345}}*-{\blacklozenge}
(-@{--}[ul],-[l])
)
,-@{--}[r]*D(3){v_3}*-{\disc}
,-@{--}[ddr]*U(3){v_5}*{\disc}(-[l],-[ur])
)
}
\end{picture}
\end{minipage}
}
\caption{The Falk arrangement $\cF_1$ and 
its associated graph}
\label{fig:falk1}
\end{figure}

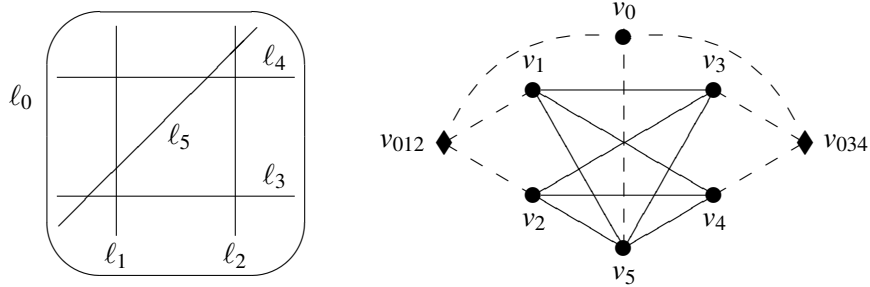
\begin{figure}%[ht]
\subfigure{%
\label{fig:f2}%
\begin{minipage}[t]{0.4\textwidth}
\setlength{\unitlength}{15pt}
\begin{picture}(5,5.8)(-3.6,-1.5)
\put(2,2){\oval(6.5,6.3)[t]}
\put(2,2){\oval(6.5,7)[b]}
\put(-1,0.5){\line(1,0){6}}
\put(-1,3.5){\line(1,0){6}}
\put(0.5,-0,5){\line(0,1){5.3}}
\put(3.5,-0,5){\line(0,1){5.3}}
\put(-0.95,-0.25){\line(1,1){5}}
\put(-1.9,3){\makebox(0,0){$\ell_0$}}
\put(0.5,-1){\makebox(0,0){$\ell_1$}}
\put(3.5,-1){\makebox(0,0){$\ell_2$}}
\put(4.5,1){\makebox(0,0){$\ell_3$}}
\put(4.5,4){\makebox(0,0){$\ell_4$}}
\put(2.1,2){\makebox(0,0){$\ell_5$}}
\end{picture}
\end{minipage}
}
\subfigure{%
\label{fig:f2graph}%
\begin{minipage}[t]{0.41\textwidth}
\setlength{\unitlength}{18pt}
\begin{picture}(5,6)(-0.3,-3.9)
\xygraph{!{0;<12mm,0mm>:<0mm,7mm>::}
!~:{@{-}|@{~}}
[]*D(3){v_0}*-{\disc}
(
-@/^-1.3pc/@{--}[ddll]*R(1.8){v_{012}}*-{\blacklozenge}
(
-@{--}[ur]
,-@{--}[dr]
)
,-@/^1.3pc/@{--}[ddrr]*L(1.8){v_{034}}*-{\blacklozenge}
(
-@{--}[ul]
,-@{--}[dl]
)
,-@{--}[dddd]*U(3){v_5}*{\disc}
(
-[uuul],-[ul],-[uuur],-[ur]
)
,[dl]*D(3){v_1}*-{\disc}
(
-[rr]
)
,[dr]*D(3){v_3}*{\disc}
(
-[ddll]
)
,[dddl]*U(3){v_2}*-{\disc}
(
-[rr]
)
,[dddr]*U(3){v_4}*{\disc}
(
-[uull]
)
)
}
\end{picture}
\end{minipage}
}
\caption{The Falk arrangement $\cF_2$ and 
its associated graph}
\label{fig:falk2}
\end{figure}

By \fullref{thm:alex poly arr}, the fundamental groups, 
$G_i=\pi_1(M(\cF_i))$, of the boundary manifolds of these 
arrangements have Alexander polynomials
\begin{align}
\label{eq:falk alex polys}
\Delta_1&=  
[(t_1{-}1)(t_2{-}1)(t_3{-}1)(t_4{-}1)(t_5{-}1)(t_{[5]}{-}1)(t_{345}{-}1)]^2
\\
\qquad\text{and}\quad\Delta_2&=
[(t_1{-}1)(t_2{-}1)(t_3{-}1)(t_4{-}1)]^2(t_5{-}1)^3(t_{[5]}{-}1)(t_{345}{-}1)
(t_{125}-1),\notag
\end{align}
where $\Delta_i=\Delta(G_i)$.  
Since these polynomials have different numbers of distinct factors, 
there is no monomial isomorphism of 
$\Lambda=\C[t_1^{\pm 1},\dots ,t_5^{\pm 1}]$ taking $\Delta_1$ 
to $\Delta_2$.  Hence, the groups $G_1$ and $G_2$ are not isomorphic, 
and the boundary manifolds $M(\cF_1)$ and $M(\cF_2)$ are not 
homotopy equivalent.  It follows that the complements of the two 
Falk arrangements are not homeomorphic---a result obtained 
previously by Jiang and Yau \cite{JY98} by invoking the 
classification of Waldhausen graph manifolds.
\end{example}

Note that the number of distinct factors in the Alexander polynomial 
$\Delta(G_2)$ above is equal to the number of vertices in the graph 
$\Gamma_{\!{\cF_2}}$, while $\Delta(G_1)$ has fewer factors than 
$|\cV(\Gamma_{\!{\cF_1}})|$.  In general, the cardinality of $\cV(\GA)$ 
is equal to $|\sD(\A)|$, the number of dense edges of $\A$.  We 
record several families of arrangements for which the Alexander 
polynomial $\Delta(G)$ is ``degenerate'', that is, the number of distinct factors 
is less than the number of dense edges.

\begin{example} 
\label{ex:degenerate}  
Let $\A$  be a line arrangement in $\CP^2$, with boundary 
manifold $M$, and let $G=\pi_1(M)$.  If $I=\{i_1,\dots,i_k\}$, 
recall that $t_I = t_{i_1}\ldots t_{i_k}$.  In particular, write 
$t_{[k]}=t_1\ldots t_k$ and $t_{[i,j]}=t_i t_{i+1} \ldots t_{j-1} t_j$.  
If $Q$ is a defining polynomial for $\A$, 
order the lines of $\A$ (starting with $0$) as indicated in $Q$.

\begin{enumerate}
\item If $Q = x_1^{n+1}-x_2^{n+1}$, then $\A$ is a pencil with 
$\abs{\sD(\A)}=n+1$ dense edges, and $G=F_{n}$ is a free group 
of rank $n$. Thus, $\Delta(G)=0$ if $n\ne 1$, and $\Delta(G)=1$ if $n=1$. 
\item If $Q = x_0(x_1^n-x_2^n)$, where $n \ge 3$, then $\A$ is a near-pencil with 
$\abs{\sD(\A)}=n+2$, while 
$\Delta(G)= (t_{[n]} - 1)^{n-2}$ has a single (distinct) factor.
\item \label{item:3}
If $Q=x_0(x_0^m-x_1^m)(x_0^n-x_2^n)$, where $m,n\ge 3$, then $\abs{\sD(\A)}=m+n+3$.
Writing $J=[m+1,m+n]$,  
$\Delta(G)$ is given by
\[
[(t_1-1)\ldots (t_m-1)(t_{[m]}-1)]^{n-1}
[(t_{m+1}-1)\ldots (t_{m+n}-1)(t_{J}-1)]^{m-1}. 
\]
\item \label{item:4}
If $Q=x_0(x_0^m-x_2^m)(x_1^n-x_2^n)$, where $m,n\ge 3$, then $\abs{\sD(\A)}=m+n+3$.  
Writing $J=[m+1,m+n]$ and $k=m+n-3$,  
$\Delta(G)$ is given by
\[
[(t_1-1)\ldots (t_m-1)(t_{[m+n]}-1)]^{n-1}[(t_{m+1}-1)\ldots 
(t_{m+n}-1)]^m (t_{J}-1)^{k}.
\]
Note that, after a change of coordinates, the Falk arrangement 
$\cF_1$ is of this form.
\end{enumerate}
\end{example}

The arrangements recorded in \fullref{ex:degenerate} \eqref{item:3} 
and \eqref{item:4} have the property that there are two $0$--dimensional 
dense edges which exhaust the lines of the arrangement.  That is, there 
are edges $F=\bigcap_{i\in I} \ell_i$ and $F'=\bigcap_{i\in I'} \ell_i$ 
so that $\A=\{\ell_i \mid i \in I \cup I'\}$.  We say $F$ and $F'$ cover $\A$.  
This condition insures that the Alexander polynomial is degenerate.

\begin{prop} 
\label{prop:degenerate}
Let $\A$ be an arrangement of $n+1$ lines in $\CP^2$ 
that is not a pencil or a near-pencil.  If $\A$ has two 
$0$--dimensional dense edges which cover $\A$, then 
the number of distinct factors in the Alexander polynomial 
of the boundary manifold of $\A$ is $\abs{\sD(\A)}-1$.  Otherwise, 
the number of distinct factors is $\abs{\sD(\A)}$.
\end{prop}

\begin{proof}
If $\A$ satisfies the hypotheses of the proposition, it is 
readily checked that, up to a coordinate change, $\A$ 
is one of the arrangements recorded in 
\fullref{ex:degenerate} \eqref{item:3} and \eqref{item:4}.  
So assume that these hypotheses do not hold.

If $\A$ has no $0$--dimensional dense edges, then $\A$ 
is a general position arrangement.  Since $\A$ is, by 
assumption, not a near-pencil, the cardinality of $\A$ 
is at least $4$, that is, $n \ge 3$.  In this instance, the 
Alexander polynomial of the boundary manifold,
\[
\Delta(G)=[(t_1-1) \ldots (t_n-1) (t_{[n]}-1)]^{n-2},
\]
has $n+1=\abs{\sD(\A)}$ factors.

Suppose $\A$ has one $0$--dimensional dense edge.  
Since $\A$ is not  a pencil or near pencil, there are at 
least two lines of $\A$ which do not contain the dense 
edge.  Write $\A=\{\ell_0,\ell_1,\dots,\ell_n\}$, where 
$\bigcap_{i=1}^k \ell_i$, $k \ge 3$, is the unique 
$0$--dimensional dense edge.  Since $\A$ has a 
single $0$--dimensional dense edge, the subarrangement 
$\{\ell_0,\ell_{k+1},\dots,\ell_n\}$ is in general position.  
By \fullref{thm:alex poly arr}, the Alexander 
polynomial of the boundary of $\A$ is
\[
\Delta(G)=\prod_{i=1}^n (t_i-1)^{m_i-2}\cdot (t_{[n]}-1)^{m_0-2} 
\cdot (t_{[k]}-1)^{k-2},
\]
and one can check that $m_i \ge 3$ for each $i$, $0 \le i \le n$.

Now consider the case where $\A$ has two $0$--dimensional 
dense edges, but they do not cover $\A$.  Either there is a 
line of $\A$ containing both dense edges, or not.  Assume 
first there is no such line.  Write $\A=\{\ell_i\}_{i=0}^n$, and 
assume without loss that the two dense edges are 
$\bigcap_{i=0}^k \ell_i$ and $\bigcap_{i=k+1}^m \ell_i$, 
where $k\ge 2$, $m-k \ge 3$, and $m<n$.  By 
\fullref{thm:alex poly arr}, 
\[
\Delta(G)=\prod_{i=1}^n (t_i-1)^{m_i-2}\cdot (t_{[n]}-1)^{m_0-2} 
\cdot (t_{[k+1,n]}-1)^{k-1}\cdot (t_{[k+1,m]}-1)^{m-k-2},
\]
and one can check that $m_i \ge 3$ for each $i$, $0 \le i \le n$.

If there is a line of $\A$ containing both $0$--dimensional dense edges, 
we can assume that $\A=\{\ell_i\}_{i=0}^n$, and the two dense edges 
are $\smash{\bigcap_{i=0}^k \ell_i}$ and $\smash{\ell_0 \cap
\bigcap_{i=k+1}^m \ell_i}$, 
where $k \ge 2$, $m-k\ge 2$, and $m < n$.  
By \fullref{thm:alex poly arr}, 
\[
\Delta(G)=\prod_{i=1}^n (t_i-1)^{m_i-2}\cdot (t_{[n]}-1)^{m_0-2} 
\cdot (t_{[k+1,n]}-1)^{k-1} \cdot (t_{[k]}t_{[m+1,n]}-1)^{m-k-1},
\]
and one can check that $m_i \ge 3$ for each $i$, $0 \le i \le n$.

Finally, suppose that $\A=\{\ell_i\}_{i=0}^n$ has at least three 
$0$--dimensional dense edges.  If $\smash{\bigcap_{i=0}^k \ell_i}$ is a 
dense edge, this assumption implies that $\smash{\bigcap_{i=k+1}^n \ell_i}$ 
cannot be a dense edge.  Consequently, the factors of the Alexander 
polynomial corresponding to $0$--dimensional dense edges are 
relatively prime, and are prime to the factor $(t_{[n]}-1)^{m_0-2}$ 
corresponding to the line $\ell_0$ of $\A$.  

To complete the 
argument, it suffices to show that $m_i \ge 2$ for each $i$, 
$0\le i\le n$.  For a  
line $\ell_i$ of $\A$, this may be 
established by choosing  
$0$--dimensional dense edges 
$F_1,F_2,F_3$ of $\A$, and considering whether $F_j$ is 
contained in $\ell_i$ or not.
\end{proof}

\section{Alexander balls} 
\label{sec:alex balls}
Let $M$ be a $3$--manifold with positive first Betti number, 
and let $G=\pi_1(M)$.  Let $H=H_1(M)/\Tors(H_1(M))$, 
and denote by $\alpha\colon G\surj H$ the projection 
onto the maximal torsion-free abelian quotient.  
Write $\rank(H)=n$, and identify $\bF[H] \cong \Lambda = 
\bF[t_1^{\pm 1},\dots,t_n^{\pm 1}]$.  Let $\phi\colon G\to \GL_k(\bF)$ 
be a linear representation, and $\Delta^{\phi}=\Delta^\phi_1(G)$ 
the corresponding twisted Alexander polynomial.  Assume 
that $\Delta^{\phi} \neq 0$, and write 
$\Delta^{\phi} =\sum c_i g_i$, 
where $0\neq c_i \in \bF$ and $g_i \in H$.  

Following McMullen \cite{Mc}, we use the twisted Alexander 
polynomial $\Delta^{\phi}$ to define a norm on 
$H^1(M;\R)=\Hom(H_1(M),\R)$.  For 
$\xi \in H^1(M;\R)$, define
\begin{equation} 
\label{eq:alex norm}
\norm{ \xi }^\phi_A := \sup_{i,j} \xi(g_i-g_j),
\end{equation}
the supremum over all $\set{g_i,g_j}$ for which 
$c_i c_j \neq 0$.  This defines a seminorm on 
$H^1(G;\R)$, the \emph{twisted Alexander norm} 
of $M$ and $\phi$.  The unit ball $\bB^\phi_A$ in 
the twisted Alexander norm is the polytope dual 
to $\cN(\Delta^{\phi})$, the Newton polytope of the 
twisted Alexander polynomial $\Delta^{\phi}$.

One also has the \emph{Thurston norm} on $H^1(M;\R)$.  
If $\Sigma$ is a compact, connected surface, let 
$\chim(\Sigma)=-\chi(\Sigma)$ if $\chi(\Sigma) \le 0$, 
and set $\chim(\Sigma)=0$ otherwise.  If $\Sigma$ is 
a surface with connected components $\Sigma_i$, set 
$\chi_{\mathunderscore}(\Sigma)=\sum\chim(\Sigma_i)$.
For $\xi \in H^1(M)$, define
\begin{equation} 
\label{eq:Thurston norm}
\norm{ \xi }_T^{\,}:= \inf\set{\chim(\Sigma)\mid \Sigma\  
\text{dual to}\ \xi},
\end{equation}
the infimum over all properly embedded oriented 
surfaces $\Sigma$.  The Thurston norm extends 
continuously to $H^1(M;\R)$.  Let $\bB_T^{\,}$ 
denote the unit ball in the Thurston norm, a polytope 
in $H^1(M;\R)$.

As shown by Friedl and Kim \cite{FK}, extending a 
result of McMullen \cite{Mc}, the twisted Alexander 
norm provides a lower bound for the Thurston norm, 
\[
\tfrac{1}{k}\norm{ \bullet }^\phi_A \le \norm{ \bullet }_T^{\,}.
\]
Consequently, the unit ball in the Thurston norm is contained 
in the unit ball in the twisted Alexander norm, 
$\bB_T^{\,} \subset \bB_A^\phi$.  For certain link 
complements, one can exhibit representations $\phi$ 
for which the twisted Alexander ball $\bB_A^\phi$ 
differs from $\bB_A$, the unit ball in the (classical) 
Alexander norm \cite{FK}, thereby distinguishing the 
Alexander and Thurston norms.  Such a distinction is 
not possible in the case where $M$ is the boundary 
manifold of a line arrangement.

\begin{thm} 
\label{thm:same alex ball}
Let $\A$ be an essential line arrangement in $\CP^2$, with boundary 
manifold $M$.  If $\phi_1$ and $\phi_2$ are complex 
representations of the  
group $G=\pi_1(M)$, 
then the twisted Alexander balls $\bB^{\phi_1}_A$ and 
$\bB^{\phi_2}_A$ are equivalent.
\end{thm}

\begin{proof}
Let $\phi\colon G\to\GL_k(\C)$ be a representation.  
We will show that the twisted Alexander ball $\bB^{\phi}_A$ 
and the (classical) Alexander ball $\bB_A$ are equivalent.  
Let $\Delta=\Delta(G)$ be the Alexander polynomial and 
$\Delta^\phi=\Delta^\phi_1(G)$ be the twisted Alexander 
polynomial associated to the representation $\phi$.  Since 
the Alexander balls $\bB_A$ and $\bB^{\phi}_A$ are the 
polytopes dual to the respective Newton polytopes of the 
Alexander polynomials, it suffices to show that $\cN(\Delta)$ 
and $\cN(\Delta^\phi)$ are equivalent.  

By \fullref{thm:alex poly arr}, the Alexander polynomials 
$\Delta$ and $\Delta^\phi$ are given by
\[
\Delta = \prod_{v \in \cV(\GA)} (t_v-1)^{m_v-2}\quad\text{and}\quad
\Delta^\phi= \prod_{v \in \cV(\GA)} \bigl[ p(A_v,t_v) \bigr]^{m_v-2},
\]
where $t_v$ is the image of a generator $z_v$ of the center of 
the vertex group $G_v$ under abelianization, and $p(A_v,t_v)$ 
is the characteristic polynomial of the automorphism $A_v=\phi(z_v)$ 
in the variable $t_v$.  Observe that only the variables $t_1,\dots,t_n$ 
appear in these Alexander polynomials.  Consequently, the Newton 
polytopes lie in $\R^n = \R^n \times \set{0} \subset H^1(M;\R)$.  
Since the Alexander polynomials factor, their Newton polytopes 
are Minkowski sums, for instance, 
\[
\cN(\Delta)=\sum_{v \in \cV(\GA)} \cN\bigl[(t_v-1)^{m_v-2}\bigr],
\]
and similarly for $\cN(\Delta^\phi)$.  

Write $d_v=m_v-2$.  If $d_v>0$ and $t_v=t_1^{q_1}\ldots t_n^{q_n}$, 
the Newton polytope $\cN\bigl[(t_v-1)^{d_v}\bigr]$ is the convex 
hull of $\b{0}=(0,\dots,0)$ and $(d_v q_1,\dots,d_v q_n)$ in 
$\R^n$, a line segment.  Thus, the Newton polytope $\cN(\Delta)$ 
is a Minkowski sum of line segments, that is, a zonotope.  As such, 
it is determined by the matrix
\begin{equation} 
\label{eq:Zmatrix}
Z=\begin{pmatrix} \b{q}_1 & \cdots & \b{q}_j\end{pmatrix},
\end{equation}
where $j$ is the number of vertices $v \in \cV(\GA)$ for which 
$d_v>0$, and
$$\b{q}_i = \begin{pmatrix} d_v q_1&\cdots
  &d_v q_n\end{pmatrix}^\top$$
if $t_v=t_1^{q_1}\ldots t_n^{q_n}$.

Now consider the Newton polytope of the twisted Alexander polynomial $\Delta^\phi$, 
\[
\cN(\Delta^\phi)=\sum_{v \in \cV(\GA)} 
\cN\bigl[p(A_v,t_v)^{d_v}\bigr].
\]  
Since the characteristic polynomial 
$p(A_v,t_v)$ is monic of degree $k$, the Newton polytope 
$ \cN\bigl[p(A_v,t_v)^{d_v}\bigr]$ is the convex hull of $\b{0}$ and 
$k\cdot (d_v q_1,\dots,d_v q_n)$ if $t_v=t_1^{q_1}\ldots t_n^{q_n}$.  
Hence, the Newton polytope $\cN(\Delta^\phi)$ is the zonotope determined 
by the matrix $k\cdot Z$, which is clearly equivalent to $\cN(\Delta)$.
\end{proof}

The Alexander and Thurston norm balls arise in the context of 
Bieri--Neumann--Strebel (BNS) invariants of the group $G=\pi_1(M)$.  
Let 
\[
\bS(G)=\bigl(H^1(G;\R)\setminus\set{\b{0}}\bigr)/\R^+,
\] 
where 
$\R^+$ acts by scalar multiplication, and view points $[\xi]$ as 
equivalence classes of homomorphisms $G\to\R$.  For $[\xi] \in \bS(G)$, 
define a submonoid $G_\xi$ of $G$ by $G_\xi=\set{g\in G \mid \xi(g)\ge 0}$.  
If $K$ is a group upon which $G$ acts, with the commutator subgroup $G'$ 
acting by inner automorphisms, the BNS invariant of $G$ and $K$ is the 
set $\varSigma_{G,K}$  of all elements $[\xi]\in\bS(G)$ for which $K$ 
is finitely generated over a finitely generated submonoid of $G_\xi$.  
The set $\varSigma_{G,K}$ is an open subset of the sphere $\bS(G)$.

Let $K=G'$, with $G$ acting by conjugation.  When $G=\pi_1(M)$, 
where $M$ is a compact, irreducible, orientable $3$--manifold, 
Bieri, Neumann, and Strebel \cite{BNS} show that the BNS invariant 
$\varSigma_{G,G'}$ is equal to the projection to $\bS(G)$ of 
the interiors of the fibered faces of the Thurston norm ball $\bB_T^{\,}$.

Assume that $H_1(M)$ is torsion-free, and consider the maximal abelian 
cover $M'$ of $M$, with fundamental group $\pi_1(M')=G'$.  The first 
homology of $M'$, $B=H_1(M') = G'/G''$, admits the structure of a 
module over $\Z[H]$, where $H=G/G'$, and is known as the 
Alexander invariant of $M$. Note that the Alexander polynomial 
$\Delta(G)=\Delta(M)$ is the order of the Alexander invariant.  
As shown by Dunfield \cite{Dun}, the BNS invariant 
$\varSigma_{G,B}$ is closely related to the Alexander polynomial.

\begin{thm} 
\label{thm:BNS & alex poly}
Let $\A$ be an essential line arrangement in $\CP^2$, with boundary manifold $M$.  
Let $G$ be the fundamental group of $M$, $B=G'/G''$ the Alexander 
invariant, and $\Delta=\ord(B)$ the Alexander polynomial.  
Then the BNS invariant $\varSigma_{G,B}$ is equal to the projection 
to $\bS(G)$ of the interiors of the top-dimensional faces of the 
Alexander ball $\bB_A$.
\end{thm}

\begin{proof}
Write $\Delta = \sum c_i g_i$, where $c_i \neq 0$ and $g_i \in H=G/G'$.  
The Newton polytope $\cN(\Delta)$ is the convex hull of the $g_i$ in 
$H_1(M;\R)$.  Call a vertex $g_i$ of $\cN(\Delta)$ a ``$\pm 1$ vertex'' 
if the corresponding coefficient $c_i$ is equal to $\pm 1$.  
For an arbitrary compact, orientable $3$--manifold $M$ whose boundary, 
if any, is a union of tori, Dunfield \cite{Dun} proves that the BNS 
invariant $\varSigma_{G,B}$ is given by the projection to 
$\bS(G)$ of the interiors of the top-dimensional faces of $\bB_A$ 
which correspond to $\pm 1$ vertices of $\cN(\Delta)$.  

If $M$ is the boundary manifold of a line arrangement 
$\A\subset \CP^2$, then, as shown in the proof of 
\fullref{thm:same alex ball}, the Newton polytope 
$\cN(\Delta)$ of the Alexander polynomial is a zonotope.  
Since the factors $(t_v-1)^{m_v-2}$ of the Alexander polynomial $\Delta$ 
have leading coefficients and constant terms equal to $\pm 1$, 
\emph{every} vertex of the associated zonotope $\cN(\Delta)$ is
a $\pm 1$ vertex.  The result follows.
\end{proof}

Let $\Delta$ be the Alexander polynomial of the boundary 
manifold of a line arrangement $\A\subset\CP^2$.  Recall 
that the Newton polytope $\cN(\Delta)$ is determined by the 
$n \times j$ integer matrix $Z$ given in \eqref{eq:Zmatrix}, 
where $\abs{\A}=n+1$ and $j$ is the number of distinct 
factors in $\Delta$.  The matrix $Z$ also determines a  
``secondary'' arrangement $\cS=\set{H_i}_{i=1}^j$ of $j$ 
hyperplanes in $\R^n$, where $H_i$ is the orthogonal 
complement of the $i$th column of $Z$.  The complement 
$\R^n \setminus\bigcup_{i=1}^j H_i$ of the real arrangement 
$\cS$ is a disjoint union of connected open sets known as 
chambers.  Let $\cham(\cS)$ be the set of chambers.  
The number of chambers may be calculated by a well known result 
of Zaslavsky \cite{Zas}.  If $P(\cS,t)$ is the Poincar\'{e} 
polynomial of (the lattice of) $\cS$, then
\[
\abs{\cham(\cS)} = P(\cS,1).
\]  
The number of chambers of the arrangement $\cS$ 
determined by the matrix $Z$ is also known to be equal to the 
number of vertices of the zonotope $\cN(\Delta)$ determined by $Z$, 
see Bj\"orner, Las Vergnas, Sturmfels, White and Ziegler \cite{BLSWZ}.
Hence, we have the following corollary to \fullref{thm:BNS & alex poly}.

\begin{cor} 
\label{cor:BNS count} 
The BNS invariant $\varSigma_{G,B}$ has $P(\cS,1)$ connected components.
\end{cor}

\begin{example} 
\label{ex:falk2}
Recall the Falk arrangements $\cF_1$ and $\cF_2$ from 
\fullref{ex:falk}.  Let $G_i$ be the fundamental 
group of the boundary manifold of $\cF_i$, 
$B_i$ the corresponding Alexander invariant, etc.  
The Alexander polynomials $\Delta_i=\Delta(G_i)$ are 
recorded in \eqref{eq:falk alex polys}.  The zonotopes 
$\cN(\Delta_1)$ and $\cN(\Delta_2)$ are determined by the matrices
\[
Z_1=\begin{pmatrix}
2 & 0 & 0 & 0 & 0 & 2 & 0\\
0 & 2 & 0 & 0 & 0 & 2 & 0\\
0 & 0 & 2 & 0 & 0 & 2 & 2\\
0 & 0 & 0 & 2 & 0 & 2 & 2\\
0 & 0 & 0 & 0 & 2 & 2 & 2
\end{pmatrix}
,\quad 
Z_2=\begin{pmatrix}
2 & 0 & 0 & 0 & 0 & 1 & 0 & 1\\
0 & 2 & 0 & 0 & 0 & 1 & 0 & 1\\
0 & 0 & 2 & 0 & 0 & 1 & 1 & 0\\
0 & 0 & 0 & 2 & 0 & 1 & 1 & 0\\
0 & 0 & 0 & 0 & 3 & 1 & 1 & 1
\end{pmatrix}.
\]
The Poincar\'e polynomials of the associated secondary 
arrangements $\cS_1$ and $\cS_2$ are
\begin{align*}
P(\cS_1,t)&=1+7t+21t^2+33t^3+27t^4+9t^5 \\
\text{and}\qquad P(\cS_2,t)&=1+8t+28t^2+51t^3+47t^4+17t^5.
\end{align*}
Consequently, the BNS invariant $\varSigma_{G_1,B_1}$ 
has $P(\cS_1,1)=98$ connected components, while 
$\varSigma_{G_2,B_2}$ has $P(\cS_2,1)=152$ 
connected components.
\end{example}

\section{Cohomology ring and holonomy Lie algebra}
\label{sect:coho}

As shown in \cite{CS06}, the cohomology ring of the boundary 
manifold  $M$ of a hyperplane arrangement has a very 
special structure:  it is the ``double" of the cohomology 
ring of the complement.  For a line arrangement, this 
structure leads to purely combinatorial descriptions of the 
skew $3$--form encapsulating $H^*(M;\Z)$, and of the 
holonomy Lie algebra of $M$.

\subsection{The doubling construction}
\label{subsec:double}

Let $R$ be a coefficient ring; we will assume either 
$R=\Z$ or $R=\bF$, a field of characteristic $0$. 
Let $A=\bigoplus_{k=0}^{m} A^k$ be a graded, 
finite-dimensional algebra over $R$.  Assume 
that $A$ is graded-commutative, of finite type (that is, each graded 
piece $A^k$ is a finitely generated $R$--module), and connected 
(that is, $A^0=R$).  Let $b_k=b_k(A)$ denote the rank of $A^k$.

Let $\bar{A} = \Hom_{R}(A,R)$ be the dual of the $R$--module 
$A$, with graded pieces $\bar{A}^{k} = \Hom_{R}(A^k,R)$.  
Then $\bar{A}$ is an $A$--bimodule, with left and right 
multiplication given by $(a\cdot f) (b)= f(ba)$ and 
$(f \cdot a) (b) =f(ab)$, respectively.  Note that, if 
$a\in A^k$ and $f\in \bar{A}^{j}$, then 
$af, fa\in \bar{A}^{j-k}$. 

Following \cite{CS06}, we define the 
{\em (graded) double} of $A$ to be the graded 
$R$--algebra $\db{A}$ with underlying $R$--module 
structure the direct sum $A\oplus \bar{A}$,  
multiplication 
\begin{equation}
\label{eq:double mult}
(a,f)\cdot (b,g) = (ab,ag+fb),
\end{equation}
for $a,b\in A$ and $f,g\in \bar{A}$, and  
grading 
\begin{equation}
\label{eq:double grading}
\db{A}^{k}=A^{k} \oplus \bar{A}^{2m-1-k}.
\end{equation}

\subsection{Poincar\'{e} duality}
\label{subsec:pd}

Let $A=\bigoplus_{k=0}^{m} A^k$ be a graded algebra as 
above. We say $A$ is a {\em Poincar\'{e} duality} algebra 
(of formal dimension $m$) if the $R$--module $A^{m}$ is 
free of rank $1$ and, for each $k$, the pairing 
$A^{k} \otimes A^{m-k} \to A^{m}$ given by multiplication 
is non-singular. In particular, each graded piece $A^{k}$ 
must be a free $R$--module. 

Given a $\PD_{m}$ algebra $A$, fix a generator $\omega$ 
for $A^{m}$.  We then have an alternating $m$--form, 
$\eta_A\colon A^1 \wedge \ldots \wedge A^1 \to R$, 
defined by 
\begin{equation}
\label{eq:eta}
a_1\ldots a_{m}= \eta_A(a_1 ,\dots, a_{m}) \cdot \omega.
\end{equation}
If $A$ is $3$--dimensional, the full multiplicative structure 
of $A$ can be recovered from the form  $\eta_A$ 
(and the generator $\omega\in A^{3}$).  

The classical example of a Poincar\'{e} duality algebra 
is the rational cohomology ring, $H^*(M;\Q)$, of an 
$m$--dimensional closed, orientable manifold $M$. 
As shown by Sullivan \cite{Su75}, any rational, alternating 
$3$--form $\eta$ can be realized as $\eta=\eta_{H^*(M;\Q)}$, 
for some $3$--manifold $M$. 

\begin{lem}
\label{lem:pddouble}
Let $A=\bigoplus_{k=0}^{m} A^k$ be a graded, graded 
commutative, connected, finite-type algebra over $R=\Z$ or $\bF$. 
Assume $A$ is a free $R$--module, and $m>1$.  If $\db{A}$ is 
the graded double of $A$,  then:
\begin{enumerate}
\item $\db{A}$ is a Poincar\'{e} duality algebra over $R$, of 
formal dimension $2m-1$.  
\item If $m> 2$, then $\eta_{\db{A}}=0$. 
\item If $m=2$, then for every $a,b,c\in A^1$ and $f,g,h\in \bar{A}^2$,
\[
\eta_{\db{A}}( (a,f), (b,g) , (c,h) ) =f(bc)+ g(ca)+h(ab). 
\] 
\end{enumerate}
\end{lem}

\begin{proof}
(1)\qua The $R$--module $\db{A}^{2m-1}=\bar{A}^0$ is isomorphic 
to $R$ via the map $f \mapsto f(1)$.  Take $\omega=\bar{1}$ 
as generator of $\db{A}^{2m-1}$.  The pairing 
$\db{A}^{k} \otimes \db{A}^{2m-1-k} 
\to \db{A}^{2m-1}$ is non-singular: its adjoint, 
\[
\db{A}^{k} \to \Hom_R( \db{A}^{2m-1-k}, \db{A}^{2m-1}), \quad
(a,f) \mapsto ((b,g) \mapsto ag+fb),
\]
is readily seen to be an isomorphism.

(2)\qua If $m>2$, then  $\db{A}^1=A^1$, and $\eta_{\db{A}}$ vanishes, 
since $A^{2m-1}=0$. 

(3)\qua If $m=2$, then $\db{A}^1=A^1 \oplus \bar{A}^2$, and the expression 
for $\eta_{\db{A}}$ follows immediately from \eqref{eq:double mult}.
\end{proof}

\subsection{The double of a $2$--dimensional algebra}
\label{subsec:2dim double}

In view of the above Lemma, the most interesting case 
is when $m=2$, so let us analyze it in a bit more detail.  
Write $A=A^0 \oplus A^1 \oplus A^2$, and fix ordered bases, 
$\{\alpha_1,\dots ,\alpha_{b_1}\}$ for $A^1$ and  
$\{\beta_1,\dots , \beta_{b_2}\}$ for $A^2$.  
The multiplication map,  
$\mu\colon A^1 \otimes A^1 \to A^2$, is then given by 
\begin{equation}
\label{eq:multiplication}
\mu(\alpha_i , \alpha_j) =
\sum_{k=1}^{b_2}\mu_{i,j,k}\,  \beta_k, 
\end{equation}
for some integer coefficients $\mu_{i,j,k}$ 
satisfying $\mu_{j,i,k}=-\mu_{i,j,k}$. 

Now consider the double 
\[
\db{A}= \db{A}^0 \oplus \db{A}^1 \oplus 
\db{A}^2  \oplus \db{A}^3 =A^0 \oplus (A^1 \oplus \bar A^2) \oplus 
(A^2 \oplus \bar A^1) \oplus \bar A^0.
\]  
Pick dual bases 
$\{\bar\alpha_j\}_{1\le j \le b_1}$  for $\bar{A}^{1}$ and 
$\{\bar\beta_k\}_{1\le k \le b_2}$  for $\bar{A}^{2}$. 
The multiplication map  
$\dbl{\mu}\colon \db{A}^1 \otimes \db{A}^1 \to \db{A}^2$ 
restricts to $\mu$ on $A^1 \otimes A^1$, vanishes on 
$\bar{A}^2\otimes \bar{A}^2$, while on $A^1\otimes \bar{A}^2$,  
it is given by
\begin{equation}
\label{eq:multi}
\dbl{\mu}(\alpha_j ,\bar \beta_k) =
\sum_{i=1}^{b_1} \mu_{i,j,k}\, \bar \alpha_i. 
\end{equation}
As a consequence, we see that the multiplication maps 
$\mu$ and $\hat{\mu}$ determine one another. 

In the chosen basis for $\db{A}^1=A^1 \oplus \bar{A}^2$, 
the form $\eta_{\db{A}}\in \bigwedge^3 \db{A}^1$ 
can be expressed~as
\begin{equation}
\label{eq:eatmu}
\eta_{\db{A}} =\sum_{1\le i< j\le b_1} \sum_{k=1}^{b_2} 
\mu_{i,j,k}  \, \alpha_i\wedge \alpha_j \wedge \bar{\beta}_k.
\end{equation}
This shows again that the multiplication map $\dbl{\mu}$ 
determines, and is determined by the $3$--form $\eta_{\db{A}}$.

\subsection{The cohomology ring of the boundary}
\label{subsec:coho ring bdry}

Now let $\A=\{\ell_i\}_{i=0}^n$ be a line arrangement in $\CP^2$, 
with complement $X$, and let $A=H^*(X;\Z)$ be the integral 
Orlik--Solomon algebra of $\A$.  As is well known, 
$A=\bigoplus_{k=0}^{2} A^k$ is torsion-free, and 
generated in degree $1$ by classes $e_1,\dots ,e_n$ dual 
to the meridians $x_1,\dots , x_n$ of the decone $\dA$. 
Choosing a suitable basis $\{f_{i,k}\mid (i,k)\in \nbc_2(\dA) \}$ 
for $A^2$, the multiplication map $\mu\colon A^1\wedge A^1 \to A^2$ 
is given on basis elements $e_i, e_j$ with $i<j$ by:
\begin{equation}
\label{eq:muarr}
\mu(e_i,e_j)=\begin{cases}
 f_{i,j} &\text{if $(i,j)\in \nbc_2(\dA)$},\\
 f_{k,j} - f_{k,i} &\text{if $\exists k$ such that 
 $(k,i), (k,j) \in \nbc_2(\dA)$},\\
 0 &\text{otherwise.}
 \end{cases}
\end{equation}
The surjectivity of $\mu$ is manifest from this formula.

For $(i,k)\in \nbc_2(\dA)$, recall that 
$I(i,k)=\set{j \mid \ell_j \supset \ell_i \cap \ell_k,\ 1\le j \le n}$.  
If $J\subset [n]$, write $e_J = \sum_{j \in J} e_j$.  
Using results from \cite{CS06} and the above discussion, 
we obtain the following.

\begin{thm}
\label{thm:coho double}
Let $\A=\{\ell_i\}_{i=0}^{n}$ be a line arrangement in $\CP^2$, 
with complement $X$ and boundary manifold $M$.  Then:
\begin{enumerate}
\item \label{c1} 
$H^*(M;\Z)$ is the double of $H^*(X;\Z)$. 
\item \label{c2} 
$H^*(M;\Z)$ is an integral Poincar\'{e} duality algebra 
of formal dimension~$3$.
\item \label{c3} 
$H^*(M;\Z)$ is generated in degree $1$ if and 
only if $\A$ is not a pencil. 
\item \label{c4} 
$H^*(M;\Z)$ determines (and is determined by) the $3$--form 
$\eta_M := \eta_{H^*(M;\Z)}$, given by
\[
\eta_M =\sum_{(i,k)\in\nbc_2(\dA)}
e_{I(i,k)} \wedge e_k \wedge \bar{f}_{i,k}.
\]
\end{enumerate}
\end{thm}

\begin{proof}
\eqref{c1} If $A=H^*(X;\Z)$, then $\db{A}=H^*(M;\Z)$, 
see \cite[Theorem~4.2]{CS06}.  

\eqref{c2} 
This follows from \fullref{lem:pddouble}, 
since $A$ is torsion-free (alternatively, use Poincar\'e duality 
for the closed, orientable $3$--manifold $M$).   

\eqref{c3} 
It is enough to show that the cup-product map 
$\dbl{\mu}\colon \dbl{A}^1\otimes \dbl{A}^1\to \dbl{A}^2$ 
is surjective if and only if $\A$ is not a pencil. 

If $\A$ is a pencil, then 
$M=\sharp^n S^1\times S^2$, and so $\dbl{\mu}=0$. 

If $\A$ is not a pencil, each line $\ell_i$ with $1\le i\le n$ 
must meet another line, say $\ell_j$, also with  $1\le j\le n$. 
Then either  $(i,j)\in \nbc_2(\dA)$, in which case 
$\hat{\mu}(e_j, \bar{f}_{i,j})=\bar{e}_i$, 
or there is an index $k\le j$ such that 
$(k,i)\in \nbc_2(\dA)$, in which case 
$\hat{\mu}(e_k, \bar{f}_{k,i})=-\bar{e}_i$.  
This shows $\bar{A}^1 \subset \Image (\hat{\mu})$. 
But we know $A^2 = \Image (\mu)$, and so 
$\hat{\mu}$ is surjective.

\eqref{c4} This follows from formulas \eqref{eq:eatmu} 
and \eqref{eq:muarr}.
\end{proof}

\begin{example} 
\label{ex:etas}
We illustrate part \eqref{c4} of the above Theorem 
with some sample computations:
\[
\eta_M=\begin{cases}
0 & \text{if $\A$ is a pencil,}\\[2pt]
(\sum_{i=1}^{n} e_i)\cdot \sum_{j=2}^{n} e_{j} \bar{f}_{1,j}
& \text{if $\A$ is a near-pencil,}\\[2pt]
\sum_{1\le i<j\le n} e_i e_j \bar{f}_{i,j}
& \text{if $\A$ is a general position arrangement.}
\end{cases}
\]
\end{example}

\begin{remark}
\label{rem:comm rel}
Let $\A$ be a line arrangement in $\CP^2$ that is not a pencil.  
Then the fundamental group $G$ of the boundary manifold $M$ 
is a commutator-relators group, $M$ is a $K(G,1)$--space, and 
$H_*(G;\Z)=H_*(M;\Z)$ is torsion-free.  In this situation, the 
cup-product structure on $H^*(G;\Z)=H^*(M;\Z)$, and hence 
the $3$--form $\eta_M$, may be computed directly from the 
commutator-relators presentation given in \fullref{prop:THEpres}, 
see, for instance, Fenn and Sjerve \cite[Theorem~2.3]{FS} and
Matei and Suciu \cite[Proposition~2.8]{MS00}.
Note, however, that the 
bases for the cohomology groups of $G$ arising in this approach 
need not, in general, coincide with those obtained from the 
realization of $H^*(G;\Z)$ as a double.
\end{remark}

\subsection{Holonomy Lie algebras}
\label{subsec:holo lie}

We now turn to a different object associated to a graded 
algebra $A$.  As before, we will assume that $A$ is 
graded-commutative, connected, and of finite-type, 
and that the ground ring $R$ is either $\Z$ or $\F$, a 
field of characteristic $0$.  Denote by $A_k$ the $R$--dual 
module $\bar{A}^k=\Hom_R(A^k,R)$; note that $A_k$ is a 
free $R$--module of rank $b_k$.

The {\em holonomy Lie algebra} 
of $A$, denoted $\h(A)$, is the quotient of the free Lie 
algebra $\Lie(A_1)$  by the ideal generated by 
the image of the comultiplication map, 
$\bar\mu \colon A_2  \to 
A_1\wedge A_1 = \Lie_2(A_1)$. 
Picking generators $x_i=\bar{\alpha}_i$ for $A_1=\bar{A}^1$, 
we obtain a finite presentation,
\begin{equation}
\label{eq:holo lie}
\h(A)=\Lie( x_1,\dots, x_{b_1} ) \Big\slash  
\Big(\sum_{1\le i<j \le b_1} \mu_{i,j,k} [x_i, x_j], \text{ for }
1\le k \le b_2\Big).  
\end{equation}
Note that $\h(A)$ inherits a natural grading from 
the free Lie algebra:  all generators $x_i$ are  in 
degree $1$, while all relations are homogeneous 
of degree $2$. 

Now let $\db{A}$ be the graded double of $A$. 
Using the description of the multiplication map 
$\hat{\mu}$ from \fullref{subsec:2dim double}, we
obtain the following presentation for $\h(\db{A})$, 
solely in terms of the multiplication map 
$\mu\colon A^1 \otimes A^1 \to A^2$, 
given by \eqref{eq:multiplication}. 

\begin{lem}
\label{lem:holo double}
The holonomy Lie algebra of $\db{A}$ is the quotient of the 
free Lie algebra on degree~$1$ generators 
$\{x_i \mid 1\le i \le b_1\}$ and $\{y_k \mid 1\le k \le b_2\}$, 
modulo the Lie ideal generated~by 
\begin{align*}
&\sum_{1\le  i<j \le b_1} \mu_{i,j,k} [x_i, x_j], 
& 1\le  k  \le b_2, \\
&\sum_{1\le j \le b_1} \sum_{1\le k \le b_2} 
\mu_{i,j,k} [x_j, y_k], & 1\le i \le b_1.
\end{align*}
\end{lem}

Note that there is a canonical projection $\h(\db{A}) \to \h(A)$, 
sending $x_i \mapsto x_i$ and $y_k\mapsto 0$. The kernel 
of this projection contains $\Lie(y_1,\dots , y_{b_2})$, 
but in general the inclusion is strict.

\subsection{The holonomy Lie algebra of the boundary}
\label{subsec:holo lie bdry}

Let $\A=\{\ell_i\}_{i=0}^n$ be a line arrangement in $\CP^2$, 
with complement $X$. As shown by Kohno \cite{K}, the 
holonomy Lie algebra of $A=H^*(X;\Z)$ has presentation 
\begin{equation}
\label{eq:holo lie arr}
\h(A)=\Lie( x_1,\dots, x_{n} ) \big\slash  
\bigg(\sum_{j\in I(i,k)}  [x_{j}, x_{k}], \text{ for }
(i,k)\in \nbc_2(\dA) \bigg).  
\end{equation}
{F}rom the preceding discussion, we find an explicit  
presentation for the holonomy Lie algebra of the 
boundary manifold of a line arrangement.  

\begin{prop}
\label{prop:holo lie bdry arr}
Let $\A=\{\ell_i\}_{i=0}^n$ be a line arrangement in $\CP^2$, 
with boundary manifold $M$.  
Then the holonomy Lie algebra of $\db{A}=H^*(M;\Z)$ is the quotient of 
the free Lie algebra on degree $1$ generators 
$\{x_i \mid  1\le i \le n\}$ and $\{ y_{(i,k)} \mid (i,k)\in \nbc_2(\dA)\}$, 
modulo the Lie ideal generated by 
$$\sum_{j\in I(i,k)} [x_{j}, x_{k}],$$
for $(i,k)\in \nbc_2(\dA)$, and
$$\sum_{ k \colon (i,k)\in \nbc_2(\dA) }  
\sum_{j\in I(i,k)}  [x_{j}, y_{(i,k)}] \ - 
\sum_{ k \colon (k,i)\in \nbc_2(\dA) }  
\sum_{j\in I(k,i)}  [x_{j}, y_{(k,i)}],$$
for $1\le i \le n$.
\end{prop}

\begin{example}
\label{ex:holo gen pos}
If $\A$ an arrangement of $n+1$ lines in general position, 
then the holonomy Lie algebra 
$\h(\db{A})$ is the quotient of the free Lie algebra on 
generators $x_i$ ($1\le i\le n$) and $y_{(i,j)}$ 
($1\le i<j\le n$), modulo the Lie ideal generated by 
\begin{align*}
&[x_i,x_j],&&1\le i< j \le n,\\
&\sum\nolimits_{j<i} [x_j, y_{(j,i)}] - \sum\nolimits_{j>i} [x_j, y_{(i,j)}],
&&1\le i \le n.
\end{align*}
\end{example}

\section{Cohomology jumping loci}
\label{sect:cjl}

In this section, we discuss the characteristic varieties 
and the resonance varieties of the boundary manifold 
of a line arrangement.  

\subsection{Characteristic varieties}
\label{subsect:char var}

Let $X$ be a space having the homotopy type of 
a connected, finite-type CW--complex.  For simplicity, we will 
assume that the fundamental group $G=\pi_{1}(X)$ has 
torsion-free abelianization $H_{1}(G)=\mathbb{Z}^n$. Consider 
the character  torus $\Hom(G,\C^*)\cong (\C^*)^{n}$.
The {\em characteristic varieties} of $X$ are the 
jumping loci for the cohomology of $X$, 
with coefficients in rank~$1$ local systems over $\C$:
\begin{equation} 
\label{eq:charvar}
V^{k}_{d}(X)=\{ \phi \in \Hom(G,\C^*) \mid 
\dim H^{k}(X; \C_{\phi})\ge d\}, 
\end{equation}
where $\C_{\phi}$ denotes the abelian group 
$\C$, with $\pi_{1}(X)$--module structure given by the 
representation $\phi\colon \pi_{1}(X)\to \C^*$.  
These loci are sub\-varieties of the algebraic torus 
$(\C^*)^{n}$; they depend only on the homotopy 
type of $X$, up to a monomial isomorphism of the 
character torus. 

For a finitely presented group $G$ (with torsion-free 
abelianization), set $V^{k}_{d}(G):=V^{k}_{d}(K(G,1))$.  
We will be only interested here in the degree $1$ characteristic 
varieties.  If $G=\pi_1(X)$ with $X$ a space as above, 
then clearly $V^{1}_{d}(G)=V^{1}_{d}(X)$. 

The varieties $V^{1}_{d}(G)$ can be computed 
algorithmically from a finite presentation of the group.  
If $G$ has generators $x_i$ and relations $r_j$, let 
$J_G=\begin{pmatrix} \partial r_i /\partial x_j\end{pmatrix}$ 
be the corresponding Jacobian matrix of Fox derivatives. 
The abelianization $J_G^{\ab}$ is the {\em Alexander matrix} 
of $G$, with entries in $\Lambda=\C[t_1^{\pm 1},\dots, t_n^{\pm 1}]$, 
the coordinate ring of $(\C^*)^n$.  Then:
\begin{equation} 
\label{eq:alex matrix}
V^{1}_{d}(G) \setminus\{1\}=V (E_d(J_G^{\ab})) \setminus\{1\}.
\end{equation}
In other words, $V^{1}_{d}(G)$ consists of all those characters 
$\phi\in \Hom(G,\C^*)\cong (\C^*)^n$ for which the evaluation 
of $J_G^{\ab}$ at $\phi$ has rank less than $n-d$ (plus, possibly, 
the identity $1$). 

\subsection{Characteristic varieties of line arrangements}
\label{subsect:char var arr}

Let $\A=\{\ell_0,\dots,\ell_n\}$ be a line arrangement in 
$\CP^2$.  The characteristic varieties of the complement 
$X$ are fairly well understood.  It follows from foundational work 
of Arapura \cite {Ar} that $V^1_d(X)$ is a union of subtori 
of the character torus $\Hom(\pi_1(X), \C^*)=(\C^*)^n$, 
possibly translated by roots of unity.  Moreover,  
components passing through $1$ admit a 
completely combinatorial description.  See \cite{CS99} and Libgober and
Yuzvinsky \cite{LY00}. 

Turning to the characteristic varieties of the boundary manifold 
$M$, we have the following complete description of 
$V^1_1(M)$.  

\begin{thm}
\label{thm:cv bdry} Let $\A$ be an essential line arrangement 
in $\CP^2$, and let $G$ be the fundamental group of 
the boundary manifold $M$.  
Then 
\[
V^1_1(G) = \bigcup_{v \in \cV(\GA),m_v \ge 3} \{t_v-1=0\}.
\]
\end{thm}

\begin{proof}
By \fullref{prop:THEpres}, the group $G$ 
admits a commutator-relators presentation, with equal 
number of generators and relations.  So the Alexander 
matrix $J_G^{\ab}$ is a square matrix, which augments to zero.   
It follows that the characteristic variety $V_1^1(G)$ is the variety defined by 
the vanishing of the codimension $1$ minors of $J_G^{\ab}$.  
The ideal $I(G)=E_1(J_G^{\ab})$ of codimension $1$ minors, 
the Alexander ideal, is given by $I(G)=\mathfrak{m}^2 \cdot (\Delta(G))$, 
where $\mathfrak{m}$ is the maximal ideal of $\Z H_1(G)$, 
see McMullen \cite{Mc}.  Consequently, 
\begin{equation}
\label{eq:cv bdry}
V^1_1(G)=\{ \Delta(G) = 0\}. 
\end{equation}
On the other hand, we know from \fullref{thm:alex poly arr} 
that the Alexander polynomial of $G$ is given by 
$\Delta(G) = \prod_{v \in \cV(\GA)} (t_v-1)^{m_v-2}$.  
The conclusion follows.
\end{proof}

By \fullref{thm:cv bdry}, $V^1_1(G)$ is the union of an arrangement
of codimension $1$ subtori in $\Hom(G,\C^*)=(\C^*)^n$, indexed by the
vertices of the graph $\GA$.  We do not have an explicit description of
the varieties $V^1_d(G)$, for $d>1$.

\subsection{Resonance varieties}
\label{subsec: res var}

Let $A$ be a graded, graded-com\-mutative, connected, 
finite-type algebra over $\C$.   
Since $a\cdot a=0$ for each $a\in A^1$, 
multiplication by $a$ defines a 
cochain complex
\begin{equation}
\label{eq:aomoto}
\xymatrixcolsep{22pt}
(A,a)\colon \:\:
\xymatrix{
0 \ar[r] &A^0 \ar[r]^{a} & A^1
\ar[r]^{a}  & A^2 \ar[r]^(.46){a}&\, \cdots }.
\end{equation}
The {\em resonance varieties} of $A$ are the 
jumping loci for the cohomology of these complexes:
\begin{equation} 
\label{eq:res var}
\RR^{k}_{d}(A)=\{ a\in A^1 \mid \dim H^k(A,a) \ge d\},
\end{equation}
for $k\ge 1$ and $1\le d \le b_k(A)$.  
The sets $\RR^{k}_{d}(A)$ are homogeneous algebraic 
subvarieties of the complex vector space $A^1=\C^{b_1}$.  

We will only be interested here in the degree $1$ resonance 
varieties, $\RR^1_{d}(A)$.  Let $S=\Sym(A_1)$ be the symmetric 
algebra on the dual of $A^1$.  If  $\{x_1,\dots ,x_{b_1}\}$ 
is the basis for $A_1$ dual to the basis 
$\{\alpha_1,\dots,\alpha_{b_1}\}$ for $A^1$, then 
$S$ becomes identified with the polynomial ring 
$\C[x_1,\dots ,x_{b_1}]$.   Also, let 
$\mu\colon A^1\otimes A^1 \to A^2$ is the 
multiplication map, given by \eqref{eq:multiplication}. 
Then, as shown by Matei and Suciu \cite{MS00} (generalizing a result 
from \cite{CS99}): 
\begin{equation}
\label{eq:res from matrix}
\RR^1_{d}(A)=V(E_{d}(\Theta)), 
\end{equation}
where  $\Theta=\Theta_A$ is the $b_1 \times b_2$ matrix of 
linear forms over $S$, with entries
\begin{equation}
\label{eq:delta matrix}
\Theta_{j,k}=\sum_{i=1}^{b_1} \mu_{i,j,k} x_i. 
\end{equation}
If $X$ is a space having the homotopy type of 
a connected, finite-type CW--complex,  define the 
resonance varieties of $X$ to be those of $A=H^*(X;\C)$. 
Similarly, if $G$ is a finitely presented group, define the 
resonance varieties of $G$ to be those of a $K(G,1)$ space. 
If $G=\pi_1(X)$, then $R^1_d(G)=R^1_d(X)$.  
Furthermore, 
if $G$ is a commutator-relators group, then 
the matrix $\Theta$ above is (equivalent to) the ``linearization" 
of the (transposed) Alexander matrix $J_G^{\ab}$,  
see \cite{MS00}.   This suggests a relationship between 
$V^1_d(G)$ and $R^1_d(G)$.   For more on this, 
see~\fullref{subsec:tcone}.

\subsection{Resonance of line arrangements}
\label{subsec:res arr}

Let $\A=\set{\ell_i}_{i=0}^n$ be an arrangement of lines in $\CP^2$,
with complement $X$.  The resonance varieties of the Orlik--Solomon
algebra $A=H^*(X;\C)$, first studied by Falk \cite{Fa97}, are by now
well understood.  It follows from \cite{CS99} and from Libgober and
Yuzvinsky \cite{LY00} that $R^1_d(A)$ is the union of linear subspaces
of $A^1=\C^n$; these subspaces (completely determined by the underlying
combinatorics) have dimension at least $2$; and intersect only at $0$.

Now let $M$ be the boundary manifold, and 
$\db{A}=H^*(M;\C)$ its cohomology ring.   
Recall that $\db{A}^1=A^1 \oplus \bar{A}^2$, with basis 
$\{\alpha_i,\bar{\beta}_k\}$, where $1\le i\le b_1=n$ and  
$1\le k \le b_2=\abs{\nbc_2(\dA)}$.  Identify the ring 
$\db{S}=\Sym(\db{A}^1)$ with the polynomial ring in 
variables $\{x_i, y_k\}$.  It follows from \eqref{eq:multi} that 
the matrix $\db{\Theta}=\Theta_{\db{A}}$ has the form 
\begin{equation}
\label{eq:bmat}
\db{\Theta}  = \begin{pmatrix} 
\Phi & \Theta \\ 
-\Theta^\top & 0 \end{pmatrix}, 
\end{equation}
where $\Phi$ is the $b_1 \times b_1$ skew-symmetric 
matrix with entries 
$\Phi_{i,j} = \sum_{k=1}^{b_2} \mu_{i,j,k} y_k$.  
Using this fact, one can derive the following information 
about the resonance varieties of $M$.  
Write $\beta = 1-b_1(A) + b_2(A)$ and $\RR_d(\Phi)= V(E_d(\Phi))$. 
 
\begin{prop}[Cohen and Suciu \cite{CS06}]
\label{prop:res var double}
The resonance varieties of the doubled algebra 
$\db{A}=H^*(M;\C)$ satisfy:
\begin{enumerate}
\item \label{rr1}
$\RR^1_d(\db{A}) = \db{A}^1$ for $d \le \beta$.
\item  \label{rr2}
$\RR^1_{d}(A) \times \bar{A}^2
\subseteq \RR^1_{d+\beta}(\db{A})$.
\item  \label{rr3}
$\RR_{d}(\Phi) \times \{0\}
\subseteq \RR^1_{d+b_2}(\db{A})$. 
\end{enumerate}
\end{prop}

This allows us to give a complete characterization of 
the resonance variety $\RR^1_1(G)$, for $G$ a boundary 
manifold group. 

\begin{cor} 
\label{cor:pen-nearpen}
Let $\A=\{\ell_0,\dots ,\ell_n\}$ be a line arrangement in $\CP^2$, 
$n\ge 2$, and $G=\pi_1(M)$. Then:
\[
\RR^1_1(G)=\begin{cases}
\C^n &\text{if $\A$ is a pencil,}\\
\C^{2(n-1)} &\text{if $\A$ is a near-pencil,}\\    
\C^{b_1+b_2} &\text{otherwise.}
\end{cases}
\]
\end{cor}

\begin{proof} 
If $\A$ is a pencil, then $G=F_n$, and so $\RR^1_1(G)=\C^{n}$. 
 
If $\A$ is a near-pencil, then $G= \Z\times \pi_1(\Sigma_{n-1})$  
and a calculation yields $\RR^1_1(G)=\C^{2(n-1)}$.

If $\A$ is neither a pencil, nor a near-pencil, then $n \ge 3$, and 
a straightforward inductive argument shows that 
$\beta\ge 1$.  Consequently, $\RR^1_1(G) =H^1(G;\C)$ 
by \fullref{prop:res var double}.
\end{proof}

\subsection{A pair of arrangements}
\label{subsec:braid prod}

\begin{figure}%
\subfigure{%
\label{fig:p23}%
\begin{minipage}[t]{0.45\textwidth}
\setlength{\unitlength}{14pt}
\begin{picture}(4,6.0)(-5,-1.6)
\put(2,2){\oval(6.5,6)[t]}
\put(2,2){\oval(6.5,7)[b]}
\put(-0.8,0.2){\line(1,0){5.6}}
\put(-0.8,2){\line(1,0){5.6}}
\put(-0.8,3.8){\line(1,0){5.6}}
\put(0.8,-0,5){\line(0,1){5}}
\put(3.2,-0,5){\line(0,1){5}}
\put(-1.9,3){\makebox(0,0){$\ell_0$}}
\put(0.8,-1){\makebox(0,0){$\ell_1$}}
\put(3.2,-1){\makebox(0,0){$\ell_2$}}
\put(4.5,-0.3){\makebox(0,0){$\ell_3$}}
\put(4.6,1.5){\makebox(0,0){$\ell_4$}}
\put(4.5,3.2){\makebox(0,0){$\ell_5$}}
\end{picture}
\end{minipage}
}
\subfigure{%
\label{fig:braidpict}%
\begin{minipage}[t]{0.45\textwidth}
\setlength{\unitlength}{14pt}
\begin{picture}(4,6.0)(-4.5,-1.6)
\put(2,2){\oval(6.5,6.3)[t]}
\put(2,2){\oval(6.5,7)[b]}
\put(-1,0.5){\line(1,0){6}}
\put(-1,3.5){\line(1,0){6}}
\put(0.5,-0,5){\line(0,1){5.3}}
\put(3.5,-0,5){\line(0,1){5.3}}
\put(-0.5,-0.5){\line(1,1){4.8}}
\put(-1.9,3){\makebox(0,0){$\ell_0$}}
\put(0.5,-1){\makebox(0,0){$\ell_1$}}
\put(3.5,-1){\makebox(0,0){$\ell_2$}}
\put(4.5,1){\makebox(0,0){$\ell_3$}}
\put(4.5,3.1){\makebox(0,0){$\ell_4$}}
\put(1.75,2.5){\makebox(0,0){$\ell_5$}}
\end{picture}
\end{minipage}
}
\caption{The product arrangement 
$\A$ and the braid arrangement $\A'$}
\label{fig:prodbraid}
\end{figure}
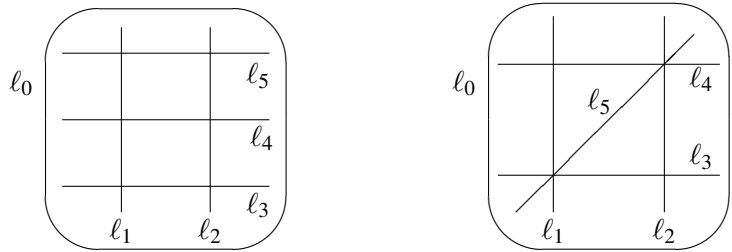

The arrangements $\A$ and $\A'$ depicted in \fullref{fig:prodbraid}
have defining polynomials
\begin{align*}
Q(\A)&=x_0(x_1+x_0)(x_1-x_0)(x_2+x_0)x_2(x_2-x_0)\\
\text{and}\qquad Q(\A')&=x_0(x_1+x_0)(x_1-x_0)(x_2+x_0)(x_2-x_0)(x_2-x_1).
\end{align*}
The respective boundary manifolds, $M$ and $M'$, 
share the same Poincar\'{e} polynomial, namely 
$P(t)=(1+t)(1+10t+t^2)$.   
Yet their cohomology rings, $\db{A}$ and $\db{A}'$, are 
not isomorphic---they 
are distinguished by their resonance varieties. Indeed, 
a computation with Macaulay~2 \cite{M2} reveals that
$$\RR^1_7(\db{A})=V(x_1,x_2,x_3,x_4,x_5,\, 
y_3y_5-y_2y_6, \, y_3y_4-y_1y_6, \, y_2y_4-y_1y_5),$$
which is a variety of dimension $4$, whereas
\begin{multline*}
\RR^1_7(\db{A}')=V(x_1,x_2,x_3,x_4,x_5,
y_2 y_4{-}y_1 y_6,\, y_2 y_5{-}y_3 y_6,
y_3 y_4{-}y_4 y_5{-}y_3 y_6+y_4 y_6,\\
 y_1 y_5{-}y_4 y_5{-}y_3 y_6+y_4 y_6,
y_1 y_3{-}y_2 y_3{-}y_4 y_5+y_1 y_6{-}y_3 y_6+y_4 y_6),
\end{multline*}
which is a variety of dimension $3$.

\section{Formality}
\label{sect:formal}
In this section, we characterize those arrangements $\A$  
for which the boundary manifold $M$ is formal, in the 
sense of Sullivan \cite{Su77}. It turns out that, with the 
exception of pencils and near-pencils, $M$ is never formal. 

\subsection{Formal spaces and $1$--formal groups}
\label{subsect:formal spaces}

Let $X$ be a space having the homotopy type of a connected, 
finite-type CW--complex.  Roughly speaking, $X$ is {\em formal}, 
if its rational homotopy type is completely determined 
by its rational cohomology ring.  More precisely, $X$ is 
formal if there is a zig-zag sequence of morphisms of 
commutative differential graded algebras connecting 
Sullivan's algebra of polynomial forms, $(A_{PL} (X,\Q),d)$,  
to $(H^*(X;\Q),0)$, and inducing isomorphisms in cohomology. 
Well known examples of formal spaces include spheres;  
simply-connected Eilenberg--Mac\,Lane spaces; 
compact, connected Lie groups and their classifying spaces;  and 
compact K\"{a}hler manifolds.  The formality property 
is preserved under wedges and products of spaces, and 
connected sums of manifolds.

A finitely presented group $G$ is said to be {\em $1$--formal}, in
the sense of Quillen \cite{Q}, if its Malcev Lie algebra (that is,
the Lie algebra of the prounipotent completion of $G$) is quadratic;
see Papadima and Suciu \cite{PS} for details.   If $X$ is a formal space,
then $G=\pi_1(X)$ is a $1$--formal group, as shown by Sullivan \cite{Su77}
and Morgan \cite{Mo}.  Complements of complex projective hypersurfaces are
not necessarily formal, see \cite{Mo}.  Nevertheless, their fundamental
groups are $1$--formal, as shown by Kohno \cite{K}.

If $X$ is the complement of a complex hyperplane arrangement, Brieskorn's
calculation of the integral cohomology ring of $X$ (see Orlik and
Terao \cite{OT1}) implies that $X$ is (rationally) formal.  However, the
analogous property of $\Z_p$--formality does not necessarily hold, due to
the presence of non-vanishing triple Massey products  in $H^*(X;\Z_p)$,
see Matei \cite{Ma}.

As mentioned above, our goal in this section is to decide, 
for a given line arrangement $\A$, whether the boundary 
manifold $M$ is formal, and whether $G=\pi_1(M)$ is 
$1$--formal. 
In our situation, Massey products in $H^*(G;\Z)$
may be computed directly from the commutator-relators
presentation given in \fullref{prop:THEpres},
using the Fox calculus approach described by Fenn and Sjerve \cite{FS}.
Yet determining whether such products vanish
is quite difficult, as Massey products are
only defined up to indeterminacy.    So we turn to other, more manageable, 
obstructions to formality.

\subsection{Associated graded Lie algebra}
\label{subsect:gr lie}

The lower central series of a group $G$ is the sequence 
of normal subgroups $\{G_k \}_{k\ge 1}$, defined 
inductively by $G_1=G$, $G_2=G'$, and $G_{k+1} =[G_k,G]$.  
It is readily seen that  the quotient groups, $G_k/G_{k+1}$,  
are abelian.  Moreover, if $G$ is finitely generated, 
so are all the LCS quotients.  
The {\em associated graded Lie algebra} of $G$ is the 
direct sum $\gr(G)=\bigoplus\nolimits_{k\ge 1} G_k/ G_{k+1}$, 
with Lie bracket  induced by the group commutator, 
and grading given by bracket length.  

If the group $G$ is finitely presented, there is another 
graded Lie algebra attached to $G$,  the (rational) 
holonomy Lie algebra, $\h(G):=\h(H^*(G;\Q))$.  In fact, 
if $X$ is any space having the homotopy type of a 
connected CW--complex with finite $2$--skeleton, 
and if $G=\pi_1(X)$, then $\h(G)=\h(H^*(X;\Q))$, see Papadima and Suciu \cite{PS}.  
Now suppose $G$ is a $1$--formal group.   Then,  
\begin{equation}
\label{eq:holo gr}
 \gr (G)\otimes \Q\cong \h(G), 
\end{equation}
as graded Lie algebras; see Quillen \cite{Q} and Sullivan \cite{Su77}.  
In particular, the respective Hilbert series must be equal. 

Returning to our situation, let $\A$ be a line arrangement 
in $\CP^2$, with boundary manifold $M$.  A finite presentation 
for the  group $G=\pi_1(M)$ is given in \fullref{prop:THEpres}.  On the other hand, 
we know that $H^*(M;\Q)=\db{A}$, 
the double of the (rational) Orlik--Solomon algebra.  Thus, 
$\h(G)=\h(\db{A})$, with presentation given in \fullref{prop:holo lie bdry arr}.  Using these explicit presentations, 
one can compute, at least in principle, the Hilbert series of  
$\gr (G)\otimes \Q$ and  $\h(G)$.

\begin{example}
\label{ex:gr holo gen pos}
Let $\A$ be an arrangement of $4$ lines in general position 
in $\CP^2$, and $M$ its boundary manifold. A presentation 
for $G=\pi_1(M)$ is given in \fullref{ex:general position}, 
while a presentation for $\h(G)$ is given in 
\fullref{ex:holo gen pos}.  Direct computation shows that 
$$\Hilb( \gr(G) \otimes \Q, t)
= 6 + 9t + 36 t^2 + 131t^3 + 528t^4 + \cdots,$$
whereas
$$\Hilb(\h(G),t) = 6 + 9t + 36 t^2 + 132 t^3 + 534 t^4 + \cdots.$$
Consequently, $G$ is not $1$--formal, and so $M$ is not formal, 
either.
\end{example}

We can use the formality test  \eqref{eq:holo gr} to 
show that several other boundary manifolds 
are not formal, but we do not know a general formula 
for the Hilbert series of the two graded Lie algebras attached 
to a boundary manifold group.  Instead, 
we turn to another formality test.

\subsection{The tangent cone formula}
\label{subsec:tcone}
Let $G$ be a finitely presented group, with $H_1(G)$ torsion-free. 
Consider the  map $\exp\colon \Hom(G,\C) \to \Hom(G,\C^*)$,  
$\exp(f)(z)=e^{f(z)}$.  Using this map, we may identify 
the tangent space at $1$ to the torus $\Hom(G,\C^*)$ 
with the vector space $\Hom(G,\C)=H^1(G,\C)$.  Under this 
identification, the exponential map takes the resonance variety  
$R^1_d(G)$ to $V^1_d(G)$. Moreover, the tangent 
cone at $1$ to $V^1_d(G)$ is contained in $R^1_d(G)$, 
see Libgober \cite{Li}.  While this inclusion is in 
general strict, equality holds under a formality assumption.

\begin{thm}[Dimca, Papadima and Suciu \cite{DPS}] 
\label{thm:tcone}
Suppose $G$ is a $1$--formal group.  Then, for each $d\ge 1$, 
the exponential map induces a complex analytic isomorphism 
between the germ at $0$ of $R^1_d(G)$ and the germ at $1$ 
of $V^1_d(G)$. Consequently,  
\begin{equation}
\label{eq:tcone}
\operatorname{TC}_{1}(V^1_d(G))=R^1_d(G).
\end{equation}
\end{thm}

In particular, this ``tangent cone formula" holds 
in the case when $X$ is the complement of a 
complex hyperplane arrangement, 
and $G$ is its fundamental group (see \cite{CS99} for 
a direct approach in this situation). 

\subsection{Formality of boundary manifolds}
\label{subsect:bdry formal}

We can now state the main result of this section, 
characterizing those line arrangements for which the 
boundary manifold is formal. 

\begin{thm}
\label{thm:nonformal}
Let $\A=\{\ell_0,\dots,\ell_n\}$ be a line arrangement in 
$\CP^2$, with boundary manifold $M$.  The following 
are equivalent:
\begin{enumerate}
\item  \label{f1} The boundary manifold $M$ is formal. 
\item  \label{f2} The group $G=\pi_1(M)$ is $1$--formal.
\item  \label{f3} The tangent cone to $V^1_1(G)$ at the 
identity is equal to $\RR^1_1(G)$. 
\item  \label{f4} $\A$ is either a pencil or a near-pencil. 
\end{enumerate}
\end{thm}

\begin{proof}
\eqref{f1} $\Rightarrow$ \eqref{f2} This follows 
from Quillen \cite{Q} and Sullivan \cite{Su77}. 

\eqref{f2} $\Rightarrow$ \eqref{f3} This follows 
from Dimca, Papadima and Suciu \cite{DPS}.

\eqref{f3} $\Rightarrow$ \eqref{f4} 
Suppose $\A$ is neither a pencil nor a near-pencil. Then  
\fullref{cor:pen-nearpen} implies that  $\RR^1_1(G) =H^1(G;\C)$. 
On the other hand, \fullref{thm:alex poly arr} implies 
that $V^1_1(G)$ is a union of codimension $1$ subtori in 
$\Hom(G,\C^*)$. Hence, the tangent cone $\operatorname{TC}_{1}(V^1_1(G))$ 
is the union of a hyperplane arrangement in $H^1(G;\C)$; thus,  
it does not equal $\RR^1_1(G)$.

\eqref{f4} $\Rightarrow$ \eqref{f1} If $\A$ is a pencil, 
then $M=\sharp^n S^1\times S^2$.  If  $\A$ is a near-pencil, 
then $M=S^1\times \Sigma_{n-1}$. In either case, $M$ 
is built out of spheres by successive product and 
connected sum operations. Thus, $M$ is formal. 
\end{proof}

\begin{rem}
The structure of the Alexander polynomial of the boundary manifold $M$
exhibited in \fullref{thm:alex poly arr} and \fullref{prop:degenerate} has
recently been used by Dimca, Papadima and Suciu \cite{DPS2} to show that
the fundamental group $G=\pi_1(M)$ is quasi-projective if and only if one
of the equivalent conditions of \fullref{thm:nonformal} holds.
\end{rem}

\begin{ack}
This research was partially supported by  National Security Agency 
grant H98230-05-1-0055 and a Louisiana State University Faculty 
Research Grant (D~Cohen), and by NSF grant DMS-0311142  
(A~Suciu). 

We thank the referee for pertinent remarks.
\end{ack}

\bibliographystyle{gtart}
\bibliography{link}

\begin{thebibliography}{}
\providecommand\bibmarginpar{\leavevmode\marginpar}
\def\urlstyle#1{{\tt #1}}

\bibitem{Ar}
\textbf{D Arapura}, \emph{Geometry of cohomology support loci for local systems
  I}, J. Algebraic Geom. 6 (1997) 563--597 \xox{MR}{1487227}

\bibitem{Ar69}
\textbf{V\,I Arnol'd}, \emph{The cohomology ring of the group of dyed braids},
  Mat. Zametki 5 (1969) 227--231 \xox{MR}{0242196}

\bibitem{BNS}
\textbf{R Bieri}, \textbf{W\,D Neumann}, \textbf{R Strebel},
  \href{http://dx.doi.org/10.1007/BF01389175} {\emph{A geometric invariant of
  discrete groups}}, Invent. Math. 90 (1987) 451--477 \xox{MR}{914846}

\bibitem{BLSWZ}
\textbf{A Bj{\"o}rner}, \textbf{M Las~Vergnas}, \textbf{B Sturmfels}, \textbf{N
  White}, \textbf{G\,M Ziegler}, \emph{Oriented matroids}, Encyclopedia of
  Mathematics and its Applications 46, Cambridge University Press, Cambridge
  (1993) \xox{MR}{1226888}

\bibitem{brown}
\textbf{K\,S Brown}, \emph{Cohomology of groups}, Graduate Texts in Mathematics
  87, Springer, New York (1982) \xox{MR}{672956}

\bibitem{CS99}
\textbf{D\,C Cohen}, \textbf{A\,I Suciu},
  \href{http://dx.doi.org/10.1017/S0305004199003576} {\emph{Characteristic
  varieties of arrangements}}, Math. Proc. Cambridge Philos. Soc. 127 (1999)
  33--53 \xox{MR}{1692519}

\bibitem{CS06}
\textbf{D\,C Cohen}, \textbf{A\,I Suciu},
  \href{http://dx.doi.org/10.1016/j.aim.2005.10.003} {\emph{Boundary manifolds
  of projective hypersurfaces}}, Adv. Math. 206 (2006) 538--566
  \xox{MR}{2263714}

\bibitem{FC}
\textbf{F\,R Cohen}, \emph{The homology of $\mathcal{C}_{n+1}$--spaces, $n>0$},
  from: ``The homology of iterated loop spaces'', Lecture Notes in Math. 533,
  Springer (1976)  207--352 \xox{MR}{0436146}

\bibitem{Dimca}
\textbf{A Dimca}, \emph{Singularities and topology of hypersurfaces},
  Universitext, Springer, New York (1992) \xox{MR}{1194180}

\bibitem{DPS2}
\textbf{A Dimca}, \textbf{S Papadima}, \textbf{A Suciu}, \emph{Alexander
  polynomials: Essential variables and multiplicities}, Int. Math. Res. Not.
  (to appear) \xox{arXiv}{0706.2499}

\bibitem{DPS}
\textbf{A Dimca}, \textbf{S Papadima}, \textbf{A Suciu}, \emph{Formality,
  Alexander invariants, and a question of Serre}  \xox{arXiv}{math.AT/0512480}

\bibitem{Dun}
\textbf{N\,M Dunfield},
  \href{http://pjm.math.berkeley.edu/pjm/2001/200-1/p03.xhtml} {\emph{Alexander
  and {T}hurston norms of fibered 3--manifolds}}, Pacific J. Math. 200 (2001)
  43--58 \xox{MR}{1863406}

\bibitem{Durfee}
\textbf{A\,H Durfee}, \href{http://dx.doi.org/10.2307/1999065}
  {\emph{Neighborhoods of algebraic sets}}, Trans. Amer. Math. Soc. 276 (1983)
  517--530 \xox{MR}{688959}

\bibitem{EN}
\textbf{D Eisenbud}, \textbf{W Neumann}, \emph{Three-dimensional link theory
  and invariants of plane curve singularities}, Annals of Mathematics Studies
  110, Princeton University Press, Princeton, NJ (1985) \xox{MR}{817982}

\bibitem{Fa}
\textbf{M Falk}, \href{http://dx.doi.org/10.1007/BF01231283} {\emph{Homotopy
  types of line arrangements}}, Invent. Math. 111 (1993) 139--150
  \xox{MR}{1193601}

\bibitem{Fa97}
\textbf{M Falk}, \href{http://dx.doi.org/10.1007/BF02558471}
  {\emph{Arrangements and cohomology}}, Ann. Comb. 1 (1997) 135--157
  \xox{MR}{1629681}

\bibitem{FS}
\textbf{R Fenn}, \textbf{D Sjerve}, \emph{Massey products and lower central
  series of free groups}, Canad. J. Math. 39 (1987) 322--337 \xox{MR}{899840}

\bibitem{FK}
\textbf{S Friedl}, \textbf{T Kim}, \emph{Twisted Alexander norms give lower
  bounds on the Thurston norm}, Trans. Amer. Math. Soc. (to appear)
  \xox{arXiv}{math.GT/0505682}

\bibitem{M2}
\textbf{D Grayson}, \textbf{M Stillman}, \emph{Macaulay~2: a software system
  for research in algebraic geometry}
\ Available at \setbox0\hbox{\makeatletter\@url
{http://www.math.uiuc.edu/Macaulay2/}}
\href{http://www.math.uiuc.edu/Macaulay2/}
{\unhbox0}

\bibitem{Hi}
\textbf{E Hironaka}, \href{http://dx.doi.org/10.1007/PL00004427}
  {\emph{Boundary manifolds of line arrangements}}, Math. Ann. 319 (2001)
  17--32 \xox{MR}{1812817}

\bibitem{Hir}
\textbf{F Hirzebruch}, \emph{The topology of normal singularities of an
  algebraic surface (after {D}. {M}umford)}, from: ``S\'eminaire Bourbaki'',
  volume 8, Exp. 250, Soc. Math. France, Paris (1995)  129--137
  \xox{MR}{1611536}

\bibitem{JY93}
\textbf{T Jiang}, \textbf{S\,S-T Yau}, \emph{Topological invariance of
  intersection lattices of arrangements in $\mathbb{C}\mathrm{P}^2$}, Bull.
  Amer. Math. Soc. $($N.S.$)$ 29 (1993) 88--93 \xox{MR}{1197426}

\bibitem{JY98}
\textbf{T Jiang}, \textbf{S\,S-T Yau},
  \href{http://www.numdam.org/item?id=ASNSP_1998_4_26_2_357_0}
  {\emph{Intersection lattices and topological structures of complements of
  arrangements in $\mathbb{C}\mathrm{P}^2$}}, Ann. Scuola Norm. Sup. Pisa Cl.
  Sci. $(4)$ 26 (1998) 357--381 \xox{MR}{1631597}

\bibitem{KL}
\textbf{P Kirk}, \textbf{C Livingston},
  \href{http://dx.doi.org/10.1016/S0040-9383(98)00039-1} {\emph{Twisted
  {A}lexander invariants, {R}eidemeister torsion, and Casson--Gordon
  invariants}}, Topology 38 (1999) 635--661 \xox{MR}{1670420}

\bibitem{KSW}
\textbf{T Kitano}, \textbf{M Suzuki}, \textbf{M Wada},
  \href{http://dx.doi.org/10.2140/agt.2005.5.1315} {\emph{Twisted {A}lexander
  polynomials and surjectivity of a group homomorphism}}, Algebr. Geom. Topol.
  5 (2005) 1315--1324 \xox{MR}{2171811}

\bibitem{K}
\textbf{T Kohno},
  \href{http://projecteuclid.org/getRecord?id=euclid.nmj/1118787354} {\emph{On
  the holonomy {L}ie algebra and the nilpotent completion of the fundamental
  group of the complement of hypersurfaces}}, Nagoya Math. J. 92 (1983) 21--37
  \xox{MR}{726138}

\bibitem{Li}
\textbf{A Libgober}, \href{http://dx.doi.org/10.1016/S0166-8641(01)00048-7}
  {\emph{First order deformations for rank one local systems with a
  non-vanishing cohomology}}, Topology Appl. 118 (2002) 159--168
  \xox{MR}{1877722}

\bibitem{LY00}
\textbf{A Libgober}, \textbf{S Yuzvinsky},
  \href{http://dx.doi.org/10.1023/A:1001826010964} {\emph{Cohomology of the
  Orlik--Solomon algebras and local systems}}, Compositio Math. 121 (2000)
  337--361 \xox{MR}{1761630}

\bibitem{Ma}
\textbf{D Matei}, \emph{Massey products of complex hypersurface complements},
  from: ``Singularity theory and its applications'', Adv. Stud. Pure Math. 43,
  Math. Soc. Japan, Tokyo (2006)  205--219 \xox{MR}{2325139}

\bibitem{MS00}
\textbf{D Matei}, \textbf{A\,I Suciu}, \emph{Cohomology rings and nilpotent
  quotients of real and complex arrangements}, from: ``Arrangements -- Tokyo
  1998'', Adv. Stud. Pure Math. 27, Kinokuniya, Tokyo (2000)  185--215
  \xox{MR}{1796900}

\bibitem{Mc}
\textbf{C\,T McMullen}, \href{http://dx.doi.org/10.1016/S0012-9593(02)01086-8}
  {\emph{The {A}lexander polynomial of a 3--manifold and the {T}hurston norm on
  cohomology}}, Ann. Sci. \'Ecole Norm. Sup. $(4)$ 35 (2002) 153--171
  \xox{MR}{1914929}

\bibitem{MT}
\textbf{G Meng}, \textbf{C\,H Taubes},
  \href{http://www.mrlonline.org/mrl/1996-003-005/1996-003-005-008.html}
  {\emph{$\underline{\mathrm{SW}}=$ {M}ilnor torsion}}, Math. Res. Lett. 3
  (1996) 661--674 \xox{MR}{1418579}

\bibitem{Mo}
\textbf{J\,W Morgan},
  \href{http://www.numdam.org/item?id=PMIHES_1978__48__137_0} {\emph{The
  algebraic topology of smooth algebraic varieties}}, Inst. Hautes \'Etudes
  Sci. Publ. Math.  (1978) 137--204 \xox{MR}{516917}

\bibitem{OT1}
\textbf{P Orlik}, \textbf{H Terao}, \emph{Arrangements of hyperplanes},
  Grundlehren der Mathematischen Wissenschaften 300, Springer, Berlin (1992)
  \xox{MR}{1217488}

\bibitem{PS}
\textbf{S Papadima}, \textbf{A\,I Suciu},
  \href{http://dx.doi.org/10.1155/S1073792804132017} {\emph{Chen {L}ie
  algebras}}, Int. Math. Res. Not.  (2004) 1057--1086 \xox{MR}{2037049}

\bibitem{Q}
\textbf{D Quillen}, \href{http://dx.doi.org/10.2307/1970725} {\emph{Rational
  homotopy theory}}, Ann. of Math. $(2)$ 90 (1969) 205--295 \xox{MR}{0258031}

\bibitem{Su75}
\textbf{D Sullivan}, \href{http://dx.doi.org/10.1016/0040-9383(75)90009-9}
  {\emph{On the intersection ring of compact three manifolds}}, Topology 14
  (1975) 275--277 \xox{MR}{0383415}

\bibitem{Su77}
\textbf{D Sullivan},
  \href{http://www.numdam.org/item?id=PMIHES_1977__47__269_0}
  {\emph{Infinitesimal computations in topology}}, Inst. Hautes \'Etudes Sci.
  Publ. Math.  (1977) 269--331 (1978) \xox{MR}{0646078}

\bibitem{Tu}
\textbf{V Turaev}, \emph{Torsions of 3--dimensional manifolds}, Progress in
  Mathematics 208, Birkh\"auser Verlag, Basel (2002) \xox{MR}{1958479}

\bibitem{Vi}
\textbf{S Vidussi},
  \href{http://pjm.math.berkeley.edu/pjm/2003/208-1/p11.xhtml} {\emph{Norms on
  the cohomology of a 3--manifold and SW theory}}, Pacific J. Math. 208 (2003)
  169--186 \xox{MR}{1979378}

\bibitem{Wa1}
\textbf{F Waldhausen}, \href{http://dx.doi.org/10.1007/BF01402956} {\emph{Eine
  Klasse von 3--dimensionalen Mannigfaltigkeiten I}}, Invent. Math. 3 (1967)
  308--333 \xox{MR}{0235576}

\bibitem{Wa2}
\textbf{F Waldhausen}, \href{http://dx.doi.org/10.1007/BF01402956} {\emph{Eine
  {K}lasse von 3--dimensionalen Mannigfaltigkeiten II}}, Invent. Math. 4 (1967)
  87--117 \xox{MR}{0235576}

\bibitem{We}
\textbf{E Westlund}, \emph{The boundary manifold of an arrangement}, PhD
  thesis, University of Wisconsin, Madison (1967)

\bibitem{Zas}
\textbf{T Zaslavsky}, \emph{Facing up to arrangements: face-count formulas for
  partitions of space by hyperplanes}, Mem. Amer. Math. Soc. 1 (1975) vii+102
  \xox{MR}{0357135}

\end{thebibliography}

\end{document}